\magnification 1200  \tolerance 10000    
\documentstyle{amsppt} 

\define\id{\operatorname{id}} 
\define\ad{\operatorname{ad}}

 \define\gr{\operatorname{gr}}

\define\ord{\operatorname{ord}}  
 
\define\End{\operatorname{End}} 
\define\Aut{\operatorname{Aut}} 
   
\define\Ss{\operatorname{\Cal S}}

\define\VGamma{\widehat{\Gamma}}  
\define\k{\Bbbk} 
\define\cbinom#1#2{\bmatrix #1 \\   #2 \endbmatrix}   

\topmatter \title \nofrills{Finite Dimensional Pointed Hopf Algebras with Abelian 
Coradical and Cartan matrices} \endtitle  
\author  {Nicol\'as Andruskiewitsch and Hans-J\"urgen Schneider} \endauthor 
\address N. Andruskiewitsch:  FAMAF. UNC. (5000) C. Universitaria. C\'ordoba.  
Argentina \endaddress 
\email {\tt andrus\@mate.uncor.edu} \endemail  
\address H.-J. Schneider: Mathematisches Seminar der Universit\"at  M\"unchen. 
Theressienstr. 39. (80333) M\"unchen. Germany \endaddress 
\email {\tt hanssch\@rz.mathematik.uni-muenchen.de}\endemail 
\thanks{(N.A.) Partially supported by CONICET, CONICOR, SeCYT (UNC)  
and   FAMAF (Rep\'ublica Argentina).} \endthanks   \endkeywords  
\leftheadtext{N. Andruskiewitsch and H.-J. Schneider} 
\rightheadtext{Pointed Hopf Algebras with Abelian Coradical}  
\date June 24, 1998\enddate
\abstract In  a previous work  \cite{AS2} we showed how to attach to a 
pointed Hopf algebra $A$ with coradical $\k\Gamma$, a braided strictly
graded Hopf algebra $R$ in the category $_{\Gamma}^{\Gamma}\Cal{YD}$ 
of Yetter-Drinfeld modules
over $\Gamma$. In this paper, we consider a further invariant of $A$,
namely the subalgebra $R'$ of $R$ generated by the space $V$ of primitive elements. 
Algebras of this kind are known since the pioneering work of Nichols.
It turns out that $R'$ is completely determined by the braiding 
$c: V\otimes V \to V \otimes V$. We denote $R' = \goth B(V)$.
We assume further that $\Gamma$ is 
finite abelian. Then $c$ is given by a matrix $(b_{ij})$ whose entries
are roots of unity; we also suppose that they have odd order. 
We introduce  for these braidings the notion of {\it braiding of Cartan type} and
we attach a generalized Cartan matrix to a braiding of Cartan type.
We prove that $\goth B(V)$ is finite dimensional if its corresponding
matrix is of finite Cartan type and give sufficient conditions for the
converse statement. As a consequence, we obtain many new families of  
pointed Hopf algebras. When $\Gamma$ is a direct sum of copies of a group
of prime order, the conditions hold and any matrix is of Cartan type.
We  apply this result to show that $R' = R$, in the case when $\Gamma$ is a group
of prime exponent. In other words, we show that a finite dimensional 
pointed Hopf algebra whose coradical is the group algebra of an abelian group
of exponent $p$ is necessarily generated by group-like and skew-primitive elements.
As a sample,  we classify all the finite dimensional coradically graded
pointed Hopf algebras whose coradical has odd prime dimension $p$. 
We also characterize coradically graded  pointed Hopf
algebras of order $p^4$.  
\endabstract \endtopmatter

\subhead \S 1. Introduction\endsubhead

Let $\k$ denote an algebraically closed field of characteristic 0.
 To motivate the results of the present paper, 
we have to recall a strategy proposed in \cite{AS2}.
 Let $A$ be a
 Hopf algebra, let $A_{0} \subseteq A_{1} \subseteq \dots$ be its coradical filtration 
and let $\gr A$ be the associated graded coalgebra. If $A_{0}$ 
is a Hopf subalgebra of $A$ (for instance, if $A$ is pointed, that is, its simple subcoalgebras are one-dimensional) then 
$\gr A = \oplus_{n\ge 0} \gr A(n)$ is a graded Hopf algebra;  $A_{0} \simeq \gr A(0)$ 
is a Hopf subalgebra and the projection $\pi: \gr A \to \gr A(0)$ with kernel
$\oplus_{n > 0} \gr A(n)$, is a Hopf algebra map and a retraction of the 
inclusion.  In this situation, a general technique due to Radford \cite{Ra} 
and explained in categorical terms by Majid \cite{Mj1} applies. Let $R$ be 
the algebra of coinvariants of $\pi$. Then $R$ is a braided Hopf algebra in 
the category $_{A_{0}}^{A_{0}}\Cal{YD}$ of Yetter-Drinfeld modules over $A_{0}$ and $A$ can be reconstructed by from $R$ and $A_{0}$. This reconstruction is called 
bosonization or biproduct.  

\medpagebreak The principle proposed in \cite{AS2} to treat questions about $A$ was  to solve the corresponding problem for $R$, to pass to $\gr A$ by bosonization and finally to lift the information to $A$. This principle was applied succesfully to the simplest possible braided Hopf algebras $R$, namely (finite dimensional) quantum linear spaces, in \cite{AS2}.   

\medpagebreak In the present paper we are concerned with braided Hopf
 algebras $R$. It turns out that a braided Hopf algebra $R$ arising from 
a coalgebra filtration as above satisfies the following conditions: 

\roster \item"(1.1)" $R = \oplus_{n\ge 0} R(n)$ is a {\it graded} braided Hopf algebra.   
\item"(1.2)" $R(0) = \k 1$ (hence the  coradical is trivial, {\it cf.} \cite{Sw, Chapter 11}).  
\item"(1.3)" $R(1) = P(R)$ (the space of primitive elements of $R$).  \endroster

 In particular, $R$ is a {\it strictly graded} coalgebra 
in the sense of \cite{Sw, Chapter 11}. 
These conditions imply that the coalgebra filtration of $R$ itself coincides with the 
filtration associated to the gradation; that is, $R$ {\it coradically graded}, 
see \cite{CM} and also \cite{AS2, Section 2}.  However, the notion of "coradically graded"
pass to $\gr A$ through bosonization but "strictly graded" not.

\medpagebreak
A first rough invariant of such $R$ is the 
dimension of $P(R)$ and we call this the {\it rank} of $R$.      

\medpagebreak It is in general not true that a braided Hopf algebra $R$ satisfying (1.1), 
(1.2) and (1.3) also verifies  

\roster \item"(1.4)" $R$ is generated as algebra over $\k$ by $R(1)$.\endroster

See Section 8. Note that the subalgebra $R'$ of $R$ generated by $R(1)$ is a Hopf 
subalgebra of $R$ and satisfies (1.1), (1.2), (1.3) and (1.4). We shall say
in general 
that a graded algebra $R= \oplus_{n\ge 0} R(n)$, with $R(0) = \k 1$, 
is {\it generated in degree one} if (1.4) holds.

\medpagebreak We are specifically interested in braided Hopf algebras 
satisfying (1.1), (1.2), (1.3) and (1.4); or in other words, which are strictly graded 
and generated in degree one.   We propose to call them Nichols algebras since they
appear for the first time in the article \cite{N}.  
Without pretending to give an exhaustive  historical analysis, let us mention that this notion was  considered by several authors  under 
various presentations; see \cite{L3}, \cite{M\"u}, \cite{Ro2}, \cite{Rz}, \cite{Sbg}, \cite{W}. 
Let us explain the reasons of our interest in this special kind of braided Hopf algebras. 
We refer to the survey article \cite{AG} for more details.

\medpagebreak First, there are alternative descriptions of Nichols algebras. 
Let $R$ be a braided Hopf 
algebra satisfying (1.1) and (1.2). Let  $V = R(1)$; it is a Yetter-Drinfeld submodule of 
$R$ contained in $P(R)$. The tensor algebra $T(V)$ is also a Yetter-Drinfeld module over 
$A_{0}$ and is in fact a braided Hopf algebra with comultiplication determined by declaring 
primitive any element of $V$. Hence the algebra morphism $\goth p: T(V) \to R$ given by the 
inclusion of $V$ on $R$  is a Hopf algebra morphism.  Now $R$ is a Nichols algebra  
if and only if 

\roster \item"(i)" $\goth p$ is surjective and 
\item"(ii)" $\ker \goth p$ is the direct sum of the kernels of the "quantum antisymmetrizers". \endroster

We recall that the "quantum antisymmetrizers" are constructed from the braiding 
$c: V\otimes V \to V\otimes V$.  However, an efficient description of the relations
does not follow easily from (ii)  and it is in general very difficult to decide when
$\goth B(V)$ is finite dimensional.
Conversely, given a Yetter-Drinfeld module $V$, 
the tensor algebra $T(V)$ divided out by the direct sum of the kernels of the 
"quantum antisymmetrizers" is a Nichols algebra with $P(R) = V$.   This line of argument also shows that a 
Nichols algebra with prescribed space of primitive elements is unique up to isomorphisms; that is, 
the notions of "Nichols algebra" and "Yetter-Drinfeld module" are naturally equivalent. We shall denote by 
$\goth B(V)$ the Nichols algebra built from the Yetter-Drinfeld module $V$.  

\medpagebreak It follows also from the description above that  
$\goth B(V)$ is not only completely determined by the Yetter-Drinfeld 
module structure of $V$,  but even more, it is completely determined by the braiding 
$c: V\otimes V \to V\otimes V$. This is a solution of the braid equation; conversely, 
given a solution of the braid equation $c: V\otimes V \to V\otimes V$, where $V$ is an 
arbitrary vector space, we can build a braided Hopf algebra 
$\goth B(V)$ in a suitable 
braided category. In this sense, Nichols algebras seem like symmetric algebras; 
however, as we shall comment below, 
they look more like enveloping algebras of Lie algebras. But there is no evident way 
(for us)  
to define a "braided Lie algebra" and, in fact, there is no need  (for our purposes) because 
all the information is in the map $c$.   

\medpagebreak Let us mention another description of Nichols algebras, 
at least when the Yetter-Drinfeld module $V$ carries an invariant non-degenerate
symmetric bilinear form. (See \cite{AG} for the general case).
Then there is a unique extension of this bilinear form to 
$T(V)$ that transforms multiplication and unit in comultiplication and counit; the kernel of 
$\goth p$ is exactly the radical of this bilinear form. Through this description,  
the positive part of the Borel-like subalgebra of a quantized enveloping algebra or of a 
Frobenius-Lusztig kernel appears as a Nichols algebra. See  \cite{L3},
\cite{M\"u}, \cite{Ro2}, \cite{Sbg}; a brief account is given in Section 3.

\medpagebreak
Finally, there is another reason supporting our interest in Nichols algebras; 
it could be possible that a 
finer understanding of them gives a way to decide whether  in characteristic 0 
a finite dimensional braided Hopf algebra $R$
 satisfying (1.1), (1.2) and (1.3) also satisfies (1.4).
See Section 8.

\medpagebreak
The main question we address in the present paper is the following: given a finite 
abelian group $\Gamma$, characterize all the  {\it  finite dimensional} 
Nichols algebras 
$R$ in $_{\Gamma}^{\Gamma}\Cal{YD}$. 
We split this question into two parts. In the first part,  we consider 
 braidings $c:V\otimes V \to V\otimes V$ arising from some finite
abelian group and discuss under which
conditions     the algebra $\goth B(V)$ is finite dimensional. 
We give several necessary and sufficient conditions; ultimately
they connect the theory of pointed of Hopf algebras, through Lusztig's work,
to Lie theory. A similar point of view is present in \cite{Ro2}, though
with different methods and results.
In the second part, it is necessary to show that a $V$ with 
$\goth B(V)$ finite dimensional actually is realizable over the fixed 
finite 
abelian group $\Gamma$. This is a problem of a somewhat different, arithmetic
nature; we completely solve it for $\Gamma = \Bbb Z/(p)$, $p$ an odd prime.

\bigpagebreak Let us assume from now on that $A_{0}$ is the group algebra of 
our fixed   finite abelian group $\Gamma$. 
Then a finite dimensional Yetter-Drinfeld module 
$V$ admits a basis $x_1$, \dots, $x_{\theta}$ such that, for some elements  
$g(1)$, \dots, $g(\theta) \in\Gamma$, $\chi(1)$, \dots, $\chi(\theta) \in\VGamma$,
 the action and coaction of $\Gamma$ are given  by 

$$ h.x_j  = \chi(j)(h) x_j, 
\qquad \delta(x_j) = g(j) \otimes x_j, \qquad j = 1,  \dots , \theta. \tag 1.5$$ 

In other words, the isomorphism class of the Yetter-Drinfeld module $V$ is 
determined (up to permutation of the index set)  by the sequences $g(1)$, 
\dots, $g(\theta) \in\Gamma$, $\chi(1)$, \dots, $\chi(\theta) \in\VGamma$.
Since we are interested in finite dimensional braided Hopf algebras $R$, we can assume that 

$$ \langle \chi(i), g(i)\rangle \neq 1, \tag 1.6$$ 
{\it cf.} \cite{AS2, Lemma 4.1}. As we said, the fundamental piece of information 
is the braiding $c$.
 Under the present hypothesis, it is given with respect to the basis
 $x_{i}\otimes x_{j}$  by 
 $$c(x_{i}\otimes x_{j}) = b_{ij} \, x_{j}\otimes x_{i}, \tag 1.7 $$   
where $(b_{ij})_{1 \le i, j\le \theta}
=(\langle \chi(j), g(i)\rangle)_{1 \le i, j\le \theta}$.

 \medpagebreak
For convenience, we shall say a braiding is represented by $(b_{ij})$ if
(1.7) holds for some basis of a fixed vector space $V$; we shall name in this case  $(b_{ij})$ the braiding  
since we reserve "matrix" for the data $(a_{ij})$ below.

 \medpagebreak  The first
main question we want to consider in this paper is: if the braiding is 
given by (1.7), when is $\goth B (V)$ finite dimensional? 
We want to make use of the theory of 
Frobenius-Lusztig kernels (see \cite{L1}, \cite{L2}, \cite{L3}). 

\definition{Definition} We shall say that a braiding given by a matrix 
$\bold b = (b_{ij})_{1\leq i,j \leq \theta}$  whose entries are roots of unity
  is {\it of Cartan type} 
if for all  $i, j$, $b_{ii}\neq 1$ 
  and   there exists  $a_{ij} \in \Bbb Z$ such that 
$$b_{ij}b_{ji} = b_{ii}^{a_{ij}}. \tag 1.8 $$ 

The integers $a_{ij}$ are completely determined
once they are chosen  in the following way: 

\roster \item"(1.9)" If $i = j$ we take $a_{ii} = 2$; 
\item"(1.10)" if $i\ne j$, we select the unique $a_{ij}$ such that 
$-\ord b_{ii} < a_{ij} \leq 0$. 
\endroster 

Then $a_{ij} = 0$ if and only if $a_{ji} = 0$,
so that $(a_{ij})$ is a generalized Cartan matrix \cite{K}.
We transfer the terminology from generalized Cartan matrices to
braidings of Cartan type:
we say that a braiding
$\bold b = (b_{ij})_{1\leq i,j \leq \theta}$ of Cartan type
is indecomposable (resp., of finite type, symmetrizable) if  $(a_{ij})$ is.
We shall also make free use of Dynkin diagrams and refer to connected Cartan
matrices for indecomposable ones.
Finally, we shall say that a Yetter-Drinfeld module $V$ is {\it of Cartan type} 
(resp., connected, ...) if the matrix  (1.7) is of Cartan type (resp., connected, ...).  
\enddefinition

\medpagebreak
We remark that there are examples of finite dimensional  braided Hopf algebras 
$\goth B (V)$ of rank 2 which are {\it not} of Cartan type, see \cite{N. pp. 1540 ff.}.  
 
\medpagebreak
The main known examples of braidings of Cartan type are given by 
$\bold b = (q^{i.j})_{i, j\in I}$, where $q$ is a root of unity
and $(I,.)$ is a Cartan datum as in \cite{L3, Chapter 1}. 
The braiding
$(b_{ij})$ is in this case symmetric. Our first reduction is to 
pass from a general braiding of Cartan type to a symmetric one; 
this is possible, at least when all the $b_{ij}$'s have odd order, 
by the twisting operation of Drinfeld. The new braiding is still of Cartan type,
and neither the elements $b_{ii}$ nor the Cartan matrix $(a_{ij})$ change.
 As we said, the algebra $\goth B(V)$ depends only on $c$ and not on the group; 
therefore, we have freedom to choose an appropiate group to perform
the twisting. See Lemma 4.1.

\remark{Remark}
There is an abuse of notation that should cause no confusion:
a symmetric braiding means that $b_{ij} = b_{ji}$ for all $i$, $j$
but neither that $a_{ij}$ is symmetric, nor that the braiding is a symmetry 
($c^2 = 1$).
\endremark
 
\medpagebreak
Even if the braiding $(b_{ij})$ is symmetric, this does not guarantee
the existence of   $q$ as above. 

\definition{Definition} Let $\bold b$ be a braiding of Cartan type with associated
Cartan matrix $(a_{ij})$ as in (1.9), (1.10). 
We say that $\bold b$ is of {\it FL-type} if there exists positive
integers $d_{1}, \dots, d_{\theta}$ such that  for all $i,j$,
\roster \item"(1.11)"   $d_{i} a_{ij} 
= d_{j}a_{ji}$ (hence $(a_{ij})$ is symmetrizable).

\item"(1.12)" There exists $q\in \k$ 
such that $b_{ij} = q^{d_{i}a_{ij}}$. \endroster  

Furthermore, we shall say that a braiding $\bold b$ is {\it locally of
FL-type} if any $2\times 2$ submatrix of $\bold b$ gives a braiding
of FL-type. \enddefinition

 \medpagebreak Let $(b_{ij})$ be a braiding of Cartan type.
Let $\Cal X$ be the set of connected components of the Dynkin diagram 
corresponding to it. For each $I \in \Cal X$, we let $N_{I}$ be the least 
common multiple of the orders of all the $b_{ii}$'s ($i\in I$), $\goth g_{I}$ 
be the Kac-Moody Lie algebra corresponding to the generalized Cartan matrix 
$(a_{ij})_{i,j \in I}$ and $\goth n_{I}$ be the Lie subalgebra of 
$\goth g_{I}$ spanned by all its  positive roots. Here is the  
main result of this article:  

\proclaim{Theorem 1.1}Let $\bold b = (b_{ij})$ be a braiding of Cartan type, 
corresponding to $c: V\otimes V \to V\otimes V$.  We also assume
that $b_{ij}$ has odd order for all $i,j$.

\roster \item If $\bold b$ is of finite type, then $\goth B(V)$ is  finite 
dimensional.  In fact $$\dim \goth B(V)
 = \prod_{I\in \Cal X}N_{I}^{\dsize \dim \goth n_{I}}.$$ 

\item  Let  $\goth B(V)$ be  finite dimensional and $\bold b$  locally of
FL-type. Let us assume that for all $i$, the order of $b_{ii}$
\endroster
\roster
\item"(a)"  is relatively 
prime to 3 whenever $a_{ij} = - 3$ for some $j$, 
\item"(b)" is different from 3, 5, 7, 11, 13, 17.

Then $\bold b$ is of finite type.\endroster \endproclaim

\medpagebreak
Let $p$ be an odd prime number. Assume that  $\bold b = (b_{ij})$ satisfies $b_{ii} \neq 1$ and
the order of $b_{ij}$ is either $p$ or 1 for all $i$ and $j$. 
Then $\bold b$ is of Cartan type.
It is also not difficult to see that it is locally of FL-type, see Lemma 4.3,
using (1.10) for $p=3$. 
We conclude from Theorem 1.1:

\proclaim{Corollary 1.2}Let $p$ be an odd prime number, 
$\Gamma$ a finite direct sum of copies of 
$\Bbb Z/(p)$ and $V$ a finite dimensional
Yetter-Drinfeld module over $\Gamma$ with braiding $\bold b$. We assume (1.6).
Then $\bold b$ is of Cartan type and
\roster \item If $\bold b$ is of finite type, then $\goth B(V)$ is  finite 
dimensional, and $\dim \goth B(V)
 = p^{M}$, $M =  \sum_{I\in \Cal X}\dsize \dim \goth n_{I}$.

\item  If $\goth B(V)$ is  finite dimensional and $p > 17$, then $\bold b$ is of finite type.
\endroster \endproclaim

\medpagebreak Now we discuss the second part of our initial question. 
Theorem 1.1 allows to construct new finite dimensional 
pointed Hopf algebras over our fixed group $\Gamma$. But we have first to
determine which matrices $\bold b$ of finite Cartan type actually appear over $\Gamma$
for some data as in (1.5).

We illustrate this in the case $\Gamma   \simeq \Bbb Z/(p)$, 
where $p$ is a prime number. 
We assume that $p$ is odd; the case $p = 2$ is considered in \cite{N, Th. 4.2.1}.

\medpagebreak
\proclaim{Theorem 1.3} Let $\Gamma \simeq \Bbb Z/(p)$, 
where $p$ is an odd prime number. 
The following  list contains  all possible
Nichols algebras of finite dimension  over $\Gamma$ and 
  finite dimensional coradically graded Hopf algebras  
with coradical isomorphic to $k\Gamma$. 

\roster \item The quantum lines and planes discussed in \cite{AS2}. By bosonization we 
get respectively Taft algebras and  book Hopf algebras {\it cf.} \cite{AS1}.

\medpagebreak \item There exists a Nichols algebra with Dynkin diagram $A_{2}$ if and only if $p = 3$ 
or $p-1$ is divisible by 3.  For $p = 3$, we obtain by bosonization from Nichols algebras of dimension 
$27$ exactly 2 non-isomorphic pointed Hopf algebras of dimension $81$ with coradical of dimension $3$.  
For $p \equiv 1 \mod 3$, we obtain by bosonization from Nichols algebras of dimension $p^{3}$ exactly $p-1$ 
non-isomorphic pointed Hopf algebras of dimension $p^{4}$ with coradical of dimension $p$.   

\medpagebreak  \item There exists a Nichols algebra with Dynkin diagram $B_{2}$ if and only if 
$p \equiv 1 \mod 4$.  For each such prime, we obtain by bosonization from Nichols algebras of 
dimension $p^{4}$ exactly $2(p-1)$ non-isomorphic pointed Hopf algebras of dimension $p^{5}$ 
with coradical of dimension $p$.   

\medpagebreak  \item There exists a Nichols algebra with Dynkin diagram $G_{2}$ if and only if 
$p \equiv 1 \mod 3$. 
For each such prime, we obtain by bosonization from Nichols algebras of dimension $p^{6}$ exactly $2(p-1)$ 
non-isomorphic pointed Hopf algebras of dimension $p^{7}$ with coradical of dimension $p$.   

\medpagebreak  \item There exist  Nichols algebras with finite Dynkin diagram of rank $\ge 3$ 
if and only if $p = 3$ and  the corresponding Dynkin diagram is $A_{2}\times A_{1}$ 
or $A_{2}\times A_{2}$.   We obtain respectively   braided Hopf algebras of dimension  
$3^{4}$ and  braided Hopf algebras of dimension  $3^{6}$ over $\Bbb Z/(3)$; hence, we  
get exactly 4 Hopf algebras of dimension $3^{5}$ and 4 Hopf algebras of dimension $3^{7}$ 
with coradical of dimension 3.  \endroster \endproclaim

These Hopf algebras are new; but the Hopf algebras of order 81 ($p = 3$, 
type $A_{2}$) was known \cite{N}.   

Theorem 1.3 opens the way to classify all finite dimensional pointed Hopf algebras
with coradical of dimension $p$ via lifting as in \cite{AS2}. 
An application of Theorem 1.3 is the determination of all possible 
 finite dimensional pointed {\it coradically graded}  Hopf algebras
of dimension $p^4$. See Theorem 7.1. Again,  the classification of  all possible 
 finite dimensional pointed Hopf algebras
of dimension $p^4$ would follow by  lifting to the filtered Hopf algebra.

 \bigpagebreak Finally, we present a striking application of Corollary 1.2. 
Theorem 1.1 would have a wider range of applications if $R$ were necessarily be
generated in degree one ({\it i.e.} $R = R'$). We see easily that this question
is related, via bosonization and a standard argument in filtered algebras- see {\it e.g.} 
\cite{AS2, Lemma 2.2}- to the following:

\proclaim{Question} Is a finite dimensional 
pointed Hopf algebra  necessarily generated by group-like and skew-primitive elements?
\endproclaim

By Theorem 1.3, the answer is yes if the coradical has prime order. 
As an application of Corollary 1.2, we prove:

\proclaim{Theorem 1.4} Let $A$ be a finite dimensional 
pointed Hopf algebra whose coradical is the group algebra of an abelian group
of prime  exponent $p >17$. Then $A$  is generated as an algebra by group-like and skew-primitive elements.
\endproclaim

We observe that this result implies that for any fixed dimension, there are only
{\it finitely many} isomorphism classes of pointed coradically graded Hopf algebras 
whose coradical is the group algebra of an abelian group
of prime  exponent $p >17$.

 \bigpagebreak The paper is organized as follows:  We discuss the necessary facts 
concerning the twisting operation (resp., Frobenius-Lusztig kernels)
in  Section 2 (resp., Section 3); we prove Theorem 1.1 in  Section 4. 
Section 5 is devoted to Nichols algebras over $\Bbb Z/(p)$ and the existence part of Theorem 1.3. 
We discuss isomorphisms between different bosonizations in Section 6; 
this concludes the proof of Theorem 1.3. Section 7 is devoted to pointed Hopf algebras
of order $p^4$.   In Section 8, we prove Theorem 1.4.

Our references for the theory of Hopf algebras are \cite{Mo}, \cite{Sw}.
Our conventions are mostly standard and were already used in \cite{AS2}.

\subsubhead Acknowledgement\endsubsubhead We are grateful to E. M\"uller for crucial conversations
about Frobenius-Lusztig kernels.
We also thank Mat\'\i as Gra\~na for 
many conversations about Nichols algebras; Jorge and Juan Guccione for enlightening 
remarks about twisting, and Sorin Dascalescu for interesting conversations.  

\bigpagebreak       

\bigpagebreak \subhead \S 2. Twisting\endsubhead

The idea of twisting was introduced by Drinfeld \cite{Dr2}, see also \cite{Dr1}, \cite{Dr3}. 
General facts are discussed in \cite{Mj2}, \cite{Mj3}; for applications to quantized 
enveloping algebras see \cite{Re}, for applications to semisimple Hopf algebras see \cite{Nk}. 
We briefly discuss the general idea referring to the literature for more details and 
work out accurately the case of our interest.  

\bigpagebreak Let $A$ be a Hopf algebra. Given an invertible $F \in A\otimes A$, we can 
consider the "twisted" comultiplication $\Delta_{F}$ given by 

$$\Delta_{F}(x) = F\Delta(x)F^{-1}.$$ 

Under certain conditions on $F$, $A_{F}$  (the same algebra $A$ but with comultiplication 
$\Delta_{F}$) is again a Hopf algebra. We are interested in the particular case when $F$ 
belongs indeed to $H\otimes H$, where $H$ is a Hopf subalgebra of $A$. We further assume 
that the inclusion $H \hookrightarrow A$ has a Hopf algebra retraction $\pi: A\to H$; so 
$A$ is the bosonization of $R$ by $H$ where $R$ is a braided Hopf algebra in the category 
$_{H}^{H}\Cal{YD}$ of Yetter-Drinfeld  modules over $H$. Then $H_{F}$ is a Hopf subalgebra 
of $A_{F}$ and $\pi$ is again a Hopf algebra map. If $H_F \simeq H$ (this happens for instance 
if $H$ is abelian), then $F$ induces an autoequivalence $V \mapsto V_{F}$ of the braided category 
$_{H}^{H}\Cal{YD}$ and $A_{F}$ is in fact the bosonization of $R_{F}$  by $H_{F}$. See \cite{Mj2}.  

\bigpagebreak We now assume that $H = \k\Gamma$ where $\Gamma$ is a finite abelian group. 
We fix  non-zero elements $y(1), \dots, y(M)$ in $\Gamma$ such that $\ord y(\ell) = E_{\ell}$,  

$$\Gamma = <y(1)> \oplus \dots \oplus <y(M)>, \quad \text{and} \quad E_{\ell}\vert E_{\ell + 1} 
\text{ for }\ell = 1, \dots, M - 1. \tag 2.1$$ 

Let $q_{i}$ be root of  unity of order $E_{i}$; we abbreviate $q = q_{M}$.  
Let $D_{i}$ be the smallest positive integer such that $q_{i} = q^{D_{i}}$.   

 \bigpagebreak Let $\gamma(i)\in \widehat{\Gamma}$ be the unique character such that 
$\langle \gamma(i),  y(j)\rangle = q_{i}^{\delta_{ij}}$, for $1\le i,j\le N$. 
Then $\ord \gamma(\ell) = E_{\ell}$ and   

$$\widehat{\Gamma} = <\gamma(1)> \oplus \dots \oplus <\gamma(M)>. $$ 

Given $g\in \Gamma$, $\chi\in \VGamma$ we shall use the notation $g_{i}$, $\chi_{j}$ meaning that  

$$ g = g_{1}y(1) + \dots + g_{M}y(M), \qquad \chi = \chi_{1}\gamma(1) + \dots + \chi_{M}\gamma(M). \tag 2.2$$  

 \bigpagebreak We identify $H$ with  the Hopf algebra $\k^{\VGamma}$ of functions on the 
group $\VGamma$; we denote by $\delta_{\tau}\in H$ the function given by $\delta_{\tau}(\zeta) 
= \delta_{\tau, \zeta}$, $\tau, \zeta\in \VGamma$. Then 

$$\delta_{\tau} = \frac 1{\vert \Gamma \vert}\sum_{g\in \Gamma} \langle \tau,  g^{-1}\rangle g.$$   

\bigpagebreak We shall describe the datum for twisting $F$ in terms of 2-cocycles, 
following \cite{Nk}, \cite{Mv}. The theory of the cocycles we are interested in goes 
back to Schur and there is an extensive literature on it; we refer to \cite{BT}.    

We recall that $\omega: \VGamma\times \VGamma \to \k^{\times}$ is a 2-cocycle if 
$\omega(\tau, 1) = \omega(1, \tau) = 1$ and  

$$\omega(\tau, \zeta)\omega(\tau\zeta, \eta) = \omega(\tau, \zeta\eta) \omega(\zeta, \eta). \tag 2.3$$ 

It is not difficult to deduce that $$ \omega(\tau\zeta, \chi)\omega(\chi, \tau\zeta)^{-1} 
\omega(\zeta, \chi)^{-1}\omega(\chi, \zeta) = \omega(\tau,\chi)\omega(\chi, \tau)^{-1}. \tag 2.4$$  

\proclaim{Lemma 2.1} Let $c_{ij}$, $1\le i < j \le M$, be integers such that $0\le c_{ij} < E_{i}$. 
Let $\omega: \VGamma\times \VGamma \to \k^{\times}$ be the map defined by   

$$ \omega(\tau, \chi) = \prod_{1 \le i < j\le M} q_{i}^{\dsize c_{ij} \tau_{j}\chi_{i}} 
= q^{\dsize\sum_{1 \le i < j\le M}  D_{i} c_{ij} \tau_{j}\chi_{i}}. \tag 2.5$$ 

Then $\omega$ is a 2-cocycle. Any 2-cocycle can be presented in this way up to a coboundary. \endproclaim 

\demo{Proof} The first statement follows by a direct computation that we leave to the reader. 
The second is not needed in the remaining of this paper, but follows from \cite{BT, Th. 4.7, p. 44}.  \qed\enddemo  

\proclaim{Lemma 2.2 \cite{Mv}, \cite{Nk}} Let $\omega: \VGamma\times \VGamma \to \k^{\times}$ 
be a 2-cocycle and let $F \in H\otimes H$ be given by  

$$F = \sum_{\tau, \zeta\in \VGamma} \omega(\tau, \zeta) \delta_{\tau} \otimes \delta_{\zeta}. $$ 

Let $A$ be a Hopf algebra containing $H$ as a Hopf subalgebra. Then $A_{F}$, the same algebra $A$ 
but with the comultiplication $\Delta_{F}$, is again a Hopf algebra. If $\widetilde \omega$ is 
another 2-cocycle and $\widetilde F$ is its corresponding twisting datum, then the Hopf algebras 
$A_{F}$ and $A_{\widetilde F}$ are isomorphic whenever $\omega$ and $\widetilde \omega$ 
differ by a coboundary. \qed\endproclaim   

\bigpagebreak We recall that a Yetter-Drinfeld  module over $H$ is a vector space $V$ provided 
with structures of left $H$-module and left $H$-comodule such that 

$$\delta(h.v) = h_{(1)} v_{(-1)} \Cal S h_{(3)} \otimes h_{(2)}.v_{(0)}.$$ 

This category is semisimple; given $g \in \Gamma$ and $\gamma \in \VGamma$, we denote by 
$V_{g}^{\gamma}$ the isotypic component of a Yetter-Drinfeld module $V$ corresponding to $g, \gamma$. 
That is, $V_{g}^{\gamma} = V_{g}\cap V^{\gamma}$, where  
$$V_{g} = \{v \in V: \delta(v) = g\otimes v\}, \qquad V^{\gamma} 
= \{v \in V: h.v = \gamma(h)v, \quad \forall h\in \Gamma\}. \tag 2.6$$ 

Conversely, a vector space $V$ provided with a direct sum decomposition 

$$V = \oplus_{g \in \Gamma, \gamma \in \VGamma} V_{g}^{\gamma}$$ 
is a Yetter-Drinfeld module with structures 
defined by (2.6). (This is a reformulation of the fact that the double of $H$ is the group algebra 
of $\Gamma\times \VGamma$).   

\bigpagebreak We now consider a braided Hopf algebra $R$ in the category $_{H}^{H}\Cal{YD}$ 
of Yetter-Drinfeld  modules over $H$ and its bosonization $A = R\# H$; then the evident maps 
$\pi: A \to H$ and $\iota: H \to A$  are Hopf algebra homomorphisms, and  
$R = \{a\in A: (\id\otimes \pi)\Delta (a) = a\otimes 1\}$. The multiplication in $A$ is 
determined by the rule $gr = \chi(g) rg$, whenever $r \in R^{\chi}$, $g \in \Gamma$; hence 

$$\delta_{\tau} r= r\delta_{\tau\chi^{-1}}, \qquad \text{if } r \in R^{\chi}.$$    

It follows from the definitions that $\pi$ and $\iota$ are Hopf algebra homomorphisms 
when replacing $A$ by its twisting $A_{F}$. Hence 

$$R_{F} := \{a\in A_{F}: (\id\otimes \pi)\Delta_{F} (a) = a\otimes 1\}$$ 

is a braided Hopf algebra in the category $_{H}^{H}\Cal{YD}$.  

\proclaim{Lemma 2.3}(i).  The linear map $\psi: R \to R_{F}$ defined by 

$$\psi (r) = \sum_{\tau\in \VGamma} \omega(\chi, \tau)^{-1}  r\# \delta_{\tau}, 
\qquad r\in R^{\chi},$$ 
is an isomorphism of  $\Gamma$-modules. If $r \in R^{\chi}$ and $s \in R^{\tau}$ then  

$$ \psi(rs) = \omega(\chi,\tau)\psi(r)\psi(s). \tag 2.7$$ 

If $r \in R^{\chi}_{g}$, then $\psi(r) \in \left(R_{F}\right)^{\chi}_{\psi(g)}$, where  

$$\aligned \psi(g) &= \sum_{\tau \in \VGamma} 
\omega(\tau, \chi)\omega(\chi,\tau)^{-1}\langle \tau,  g\rangle \delta_{\tau} 
\\ &= \sum_{1\le j \le M}   \left(\sum_{1\le i < j} \frac{D _{i}}{D_{j}} c_{ij} \chi_{i}  
+ g_{j}- \sum_{j < h\le M} c_{jh} \chi_{h}\right) y(j). \endaligned\tag 2.8$$  
Here we have used notations (2.2) and (2.5). 
(The first equality of (2.8) is in $\k\Gamma$, the second in $\Gamma$).   

\medpagebreak (ii). If $R$ is a graded braided Hopf algebra, then $R_{F}$ also is
 and $\psi$ is a graded map. If $R$ is a coradically graded braided Hopf algebra 
(resp. a Nichols algebra), then $R_{F}$ also is. \endproclaim  

\demo{Proof}  We first check that $\psi$ is well-defined: 

$$\align (\id\otimes \pi)\Delta_{F} (\psi(r)) =& F(r\otimes 1) 
\sum_{\tau, \sigma\in \VGamma} \omega(\chi, \sigma\tau)^{-1}   
\delta_{\tau} \otimes \delta_{\sigma}  F^{-1}\\ &= \sum_{\mu, \eta\in \VGamma} 
\omega(\mu, \eta) \delta_{\mu} \otimes \delta_{\eta} (r\otimes 1) 
\sum_{\tau, \sigma\in \VGamma} \omega(\chi, \tau\sigma)^{-1} 
\omega(\tau, \sigma)^{-1} \delta_{\tau} \otimes \delta_{\sigma} 
\\ &= (r\otimes 1) \sum_{\tau, \sigma\in \VGamma} \omega(\tau\chi, \sigma) 
\omega(\chi, \tau\sigma)^{-1} \omega(\tau, \sigma)^{-1} \delta_{\tau} \otimes \delta_{\sigma} 
\\ &= (r\otimes 1) \sum_{\tau, \sigma\in \VGamma}\omega(\chi, \tau)^{-1} \delta_{\tau} 
\otimes \delta_{\sigma} \\ &= \psi(r)\otimes 1. \endalign$$ 

It is easy to see that $\psi$ preserves the action of $\Gamma$. It follows that 
$\psi$ is invertible with inverse  $\psi^{-1} (r) = 
\sum_{\tau\in \VGamma} \omega(\chi, \tau) r\# \delta_{\tau}$, for $r\in \left(R_{F}\right)^{\chi}$. 

Now we prove (2.7): 

$$\align \psi(r)\psi(s) &= \sum_{\sigma, \eta\in \VGamma} \omega(\chi, \sigma)^{-1} 
\omega(\tau, \eta)^{-1} r\delta_{\sigma} s\delta_{\eta} = \sum_{\eta\in \VGamma} 
\omega(\chi, \tau\eta)^{-1} \omega(\tau, \eta)^{-1} rs\delta_{\eta} 
\\ &  =  \omega(\chi, \tau)^{-1} \sum_{\eta\in \VGamma}\omega(\chi\tau, \eta)^{-1} rs\delta_{\eta} 
= \omega(\chi, \tau)^{-1} \psi(rs). \endalign $$ 

We still assume that  $r \in R^{\chi}_{g}$. Then the coaction on $\psi(r)$ is given by 

$$ \align \delta_{F} (\psi(r)) & = (\pi\otimes \id)\Delta_{F} (\psi(r)) 
\\ & = \left(\sum_{\mu, \eta} \omega(\mu, \eta) \delta_{\mu}\otimes \delta_{\eta}\right) 
\left(g\otimes r\right) \Delta\left(\sum_{\tau} 
\omega(\chi,\tau)^{-1} \delta_{\tau} \right) F^{-1} 
\\ & = \sum_{\tau, \zeta} \omega(\tau, \zeta\chi) \omega(\chi, \tau\zeta)^{-1} 
\omega(\tau, \zeta)^{-1} g\delta_{\tau} \otimes r\delta_{\zeta}  
\\ &  = \sum_{\tau, \zeta}  \omega(\chi, \tau\zeta)^{-1} 
\omega(\zeta, \chi)^{-1} \omega(\tau\zeta, \chi) g\delta_{\tau} \otimes r\delta_{\zeta}  
\\ &  = \sum_{\tau, \zeta}  \omega(\chi, \zeta)^{-1} \omega(\tau, \chi) \omega(\chi, \tau)^{-1} 
g\delta_{\tau} \otimes r\delta_{\zeta}  \\ &  = \sum_{\tau}   \omega(\tau, \chi) 
\omega(\chi, \tau)^{-1} \langle \tau, g \rangle\delta_{\tau} \otimes\psi(r), \endalign $$ 
where we have used repeatedly the cocycle condition and (2.4). 

From the coassociativity of $\Delta_{F}$ we conclude that 
$$\psi(g) = \sum_{\tau}   \omega(\tau, \chi) \omega(\chi, \tau)^{-1} \langle \tau, g \rangle\delta_{\tau}$$ 
is a group-like element of $H$, {\it i.e.} an element of $\Gamma$. 
The last equality in (2.8) follows by evaluating at $\gamma(j)$.   

\medpagebreak Let us suppose now that $R = \oplus_{n\ge 0} R(n)$ is a {\it graded} braided Hopf 
algebra. Then $A = \oplus_{n\ge 0} A(n)$, where $A(n) = R(n) \# H$, is  a graded Hopf algebra. 
Since $F\in H\otimes H = A(0)\otimes A(0)$, we see that $A_{F}$ is also a graded Hopf algebra 
with respect to the same grading and $\pi$, $\iota$ are graded maps. Hence $R_{F}$ inherits 
the grading of $A$ and is a graded braided Hopf algebra. It is now evident that $\psi$ is 
a graded map, $R$ and $R_{F}$ being both  graded subspaces of $A$.  

\medpagebreak  Let us further suppose that $R$ is coradically graded; then $A$ also is. 
Since $F\in H\otimes H = A(0)\otimes A(0)$, $A_{F}$ has the same coradical filtration 
than $A$; hence it is also coradically graded. 
Therefore $R_{F}$ is coradically graded. See \cite{AS2, Section 2}. In particular, 
$\psi(P(R)) = P(R_{F})$.    If in addition $R$ 
is generated by $P(R)$, {\it i.e.} is a Nichols algebra, then we conclude from (2.7) that $R_{F}$ 
is generated by $P(R_{F})$, 
thus $R_{F}$ is also a Nichols algebra. \qed\enddemo  

\bigpagebreak Let now $R$ be a coradically graded braided Hopf algebra, 
$x_{1}, \dots, x_{\theta}$ a basis of $P(R)$ such that for some 
$g(1), \dots, g(\theta)\in \Gamma$, $\chi(1), \dots, \chi(\theta)\in\VGamma$ (1.5) holds. 
Then the braiding (1.7) is given by  

$$ \aligned (b_{ij})_{1 \le i, j\le \theta} &=  
(\langle \chi(i), g(j)\rangle)_{1 \le i, j\le \theta} \\ &=  
\left(\prod_{1\le h \le M}q_{h}^{\chi(i)_{h} g(j)_{h}}\right)_{1 \le i, j\le \theta} 
= (q^{\alpha_{ij}}), \\ \text{where }\quad \alpha_{ij}
&= \sum_{1\le h \le M} D_{h} \chi(i)_{h} g(j)_{h}. \endaligned \tag 2.9 $$ 

We keep the same meaning for $F$ as above. By Lemma 2.3, $R_{F}$ is a coradically
 graded braided Hopf algebra and $\psi(x_{1}), \dots, \psi(x_{\theta})$ is a basis 
of $P(R_{F})$ such that (1.5) holds for $\psi(g(1)), \dots, \psi(g(\theta))\in \Gamma$, 
$\chi(1), \dots, \chi(\theta)\in\VGamma$. By (2.8), the braiding (1.7) for $R_{F}$ is given by 

$$ \aligned (b^{F}_{ij})_{1 \le i, j\le \theta} &=  (\langle \chi(i), 
\psi(g(j))\rangle)_{1 \le i, j\le \theta} =   (q^{\alpha^{F}_{ij}}), \\ 
\\ \text{where }\quad \alpha^{F}_{ij}&= \sum_{1\le t < h \le M} D_{t} c_{th} 
\left( \chi(i)_{h} \chi(j)_{t} - \chi(i)_{t} \chi(j)_{h}\right)  + \alpha_{ij}. 
\endaligned \tag 2.10 $$

\remark{Remark 2.4} The "symmetrization" of the matrix $\alpha_{ij}$ 
remains unchanged under twisting: 

$$ \aligned \alpha^{F}_{ij} + \alpha^{F}_{ji} &= 
\sum_{1\le t < h \le M} D_{t} c_{th} \left( \chi(i)_{h} \chi(j)_{t} 
- \chi(i)_{t} \chi(j)_{h} \right. 
\\  &\qquad \left. + \chi(j)_{h} \chi(i)_{t} - \chi(j)_{t} \chi(i)_{h}\right) 
+ \alpha_{ij} + \alpha_{ji} = \alpha_{ij} + \alpha_{ji}. 
\endaligned \tag 2.11$$ 

Hence, if $R$ is of Cartan type, then $R_{F}$ also is. 
Also, $ \alpha^{F}_{ii} =  \alpha_{ii}$ and the diagonal
elements $b_{ii}$ do not change after twisting.

In particular, if $R$ is a quantum linear space, {\it i.e.}  
$\alpha_{ij} + \alpha_{ji} = 0 \mod E_{M}$ whenever $i\neq j$, 
then $R_{F}$ also is a quantum linear space.   \endremark

 \bigpagebreak \subhead \S 3. Frobenius-Lusztig kernels\endsubhead

Let $(I, \cdot)$ be a Cartan datum as in \cite{L3, Chapter 1} and let $(Y, X, \dots)$ be a root datum
of type $(I, \cdot)$. Let $a_{ij} = \langle i, j' \rangle$,
$\dsize d_{i} = \frac{i\cdot i}2$. Let also
$q \neq 1$ be a root  of 1 of odd order $N$. We consider in this Section a vector space $V$
with braiding $\bold b = (b_{ij})$
where $b_{ij} = q^{i\cdot j}$.  
We also assume that
\roster \item"(3.1)" $N_{i} > -a_{ij}$, where $N_i$ is the order of $b_{ii}$, for all $i,j$.
\endroster
Note that $N_{i}$ has the meaning with respect to $N$ as in \cite{L3, 2.2.4}.

\medpagebreak
Now $\bold b$ is of Cartan type and by (3.1), it is of FL-type and its associated Cartan
matrix is $(a_{ij})$, with the $d_{i}$'s as in the preceding sentence. Conversely,
any braiding of FL-type arises in this way.

Our aim is 
to sketch the main ideas of the proof of the following Theorem, which follows from deep results
of Lusztig and in addition from the work of M\"uller \cite{M\"u} and Rosso \cite{Ro1}.

\proclaim{Theorem 3.1} \roster \item Assume that $(I, \cdot )$ is of finite type. Then $\goth B(V)$
is isomorphic to the positive part of the Frobenius-Lusztig kernel $\goth u$, {\it cf.}
\cite{L3}. In particular, $\dim \goth B(V) = N^M$, where $M$ is the dimension of the nilpotent subalgebra
of the Lie algebra corresponding to $(I, \cdot )$ spanned by all the positive root vectors.

\item  Assume that $(I, \cdot )$ is not of finite type. Then $\goth B(V)$ is infinite dimensional. \endroster

\endproclaim

\demo{Proof} Part (1) follows from the characterization in \cite{M\"u, Section 2} or \cite{Ro1}; the
claim on the dimension follows from \cite{L1}, \cite{L2}.

We sketch now the proof of Part (2). We keep the notation from \cite{L3}. 
It is in principle not always true that the positive part $\goth u^{+}$ of the Frobenius-Lusztig kernel
is isomorphic to $\goth B(V)$. We have however an epimorphism $\goth u^{+}\to \goth B(V)$.

\medpagebreak
Since $(I, \cdot )$ is not of finite type, its Weyl group is not finite. Hence, for any $M >0$
there exists an element $w \in W$ of length $M$, say with reduced expression $w = s_{i_{1}} \dots s_{i_{M}}$.
We shall show that $\goth B(V)$ contains $2^M$ linearly independent elements; this concludes the proof of (2).

We work first in the transcendental case.
We use now specifically the notation from \cite{L3, 38.2.2}. The sequence
$\bold h = (i_{1}, \dots, i_{M})$ is admissible because of \cite{L3, Proposition 40.2.1 (a)}. 
It is clear that the element $x = 1$ is adapted to $(\bold h, 0)$.
Let $c$, $c'$ be as \cite{L3, Proposition 38.2.3} and suppose that their entries are either 0 or 1.
Then $L(\bold h, c, 0, 1)$ is orthogonal (resp, not orthogonal) to  $L(\bold h, c', 0, 1)$
if $c$ is different (resp., equal) to $c'$. But we can pass from the transcendental case to our case
($q$ a root of unity) thanks to  the results in \cite{M\"u, Section 2}, notably Lemma 2.2,
Proposition 2.3 and Theorem 2.11 (a). \qed \enddemo

\bigpagebreak \subhead \S 4. Proof of Theorem 1.1\endsubhead

\subsubhead 4.1 Sketch of the proof\endsubsubhead 

Let us first outline the proof of Theorem 1.1.  Let $V$ be a Yetter-Drinfeld module
over our fixed group $\Gamma$ as in (1.5) and let $R = \goth B(V)$.
We begin by passing from a general braiding (1.7) 
to a symmetric one. For this, it is natural to try to apply a Drinfeld twist to the Hopf 
algebra $\k \Gamma \# R$ (compare with \cite{Ro2}). 
Unfortunately  this is not always possible. 
It fails already for quantum linear planes over cyclic groups, since no twist is possible 
in the cyclic group case. However, we can overcome this difficulty because the algebra 
and coalgebra structures of $R$ do not depend on $\Gamma$ but only on the braiding (1.7).  
That is, let $\Upsilon$ be another finite abelian group and let $h(1)$, \dots, 
$h(\theta) \in\Upsilon$, $\eta(1)$, \dots, $\eta(\theta) \in\widehat{\Upsilon}$ 
be sequences satisfying 

$$\langle \eta(i), h(j)\rangle = \langle \chi(i), g(j)\rangle, \quad \forall  i,j.$$ 

Let $T$ be the Nichols algebra over $\Upsilon$ such that $T(1) = V$ with Yetter-Drinfeld module 
structure given by (1.5) but with $h(i)$'s and $\eta(j)$'s instead of $g(i)$'s 
and $\chi(j)$'s. Then $T$ is isomorphic to $R$ both as an algebra and coalgebra. 
This follows at once from the description of $T$ and $R$ in terms of 
"quantum antisymmetrizers" mentioned in the Introduction; see for instance 
\cite{Sbg}, \cite{Ro2}, \cite{AG}.  So we change the group,  twist and  
assume that the braiding (1.7) is symmetric.   

\medpagebreak
We then discuss briefly the relations between the notions of Cartan type,
FL-type and locally of FL-type. We state the last one in arithmetical terms
and describe several instances where it holds. One of them is the following:

Let  $\bold b = (b_{ij})$ be a  connected
braiding of Cartan type whose entries have odd order. We say
that $\bold b$ satisfies the {\it relative primeness condition}
if $a_{ij}$ is 0 or relatively prime to the order of $b_{ii}$ for all $i,j$.

We prove then that a braiding of symmetrizable Cartan type satisfying the relative 
primeness condition is of FL-type. By the results evoked in Section 3, this is enough
for part (1) of Theorem 1.1. 

For the converse part (2) of Theorem 1.1, we prove that if $\goth B(V)$
is finite dimensional and $\bold b$  is locally of FL-type
then it is symmetrizable. 
By the results of Section 3 again, this concludes the proof of the theorem.

\bigpagebreak 
\subsubhead 4.2 Reduction to the symmetric connected case\endsubsubhead

Let $(b_{ij})_{1\le i,j \le \theta}$ be a matrix with entries in $\k^{\times}$.  
Let $V$ be a vector space with a basis $x_{1}, \dots, x_{\theta}$ and let 
$c: V\otimes V \to V\otimes V$ be the linear map given by 

$$c(x_{i}\otimes x_{j}) = b_{ij} x_{j}\otimes x_{i}, 
\quad 1\le i,j \le \theta. \tag 4.1$$

It is easy to verify that $c$ satisfies the Braid relation: 

$$(c\otimes \id) (\id \otimes c) (c\otimes \id) 
= (\id \otimes c) (c\otimes \id) (\id \otimes c).$$ 

Therefore, it induces representations of the Braid group $\Bbb B_{n}$ 
on $V^{\otimes n}$ for any $n\ge 2$ and this in turn allows to define  
"braided antisymmetrizers" which are linear maps $\goth S_{n}: V^{\otimes n} 
\to V^{\otimes n}$. Let $R = \oplus_{n\ge 0} V^{\otimes n} / \ker \goth S_{n}$; 
it is a graded  algebra. See for instance \cite{Sbg}, \cite{AG}. Using the 
methods of \cite{Mj3, Th. 10.2.1} it is possible to prove that $R$ is a braided Hopf algebra. 
In our case this will follow from Lemma 4.1 below.  

\bigpagebreak If $\Gamma$ is an abelian group, $g(1)$, \dots, $g(\theta) \in\Gamma$, 
$\chi(1)$, \dots, $\chi(\theta) \in\widehat{\Gamma}$ are  sequences such that 
$b_{ij} = \langle \chi(j),  g(i)\rangle$ and we consider $V$ as a Yetter-Drinfeld 
module over $\Gamma$  via (1.5), then $c$ is the "commutativity isomorphism" in 
the braided category $_{\Gamma}^{\Gamma}\Cal{YD}$. Moreover, $R$ has a braided 
Hopf algebra structure and is in fact the Nichols algebra associated to $V$.    

\proclaim{Lemma 4.1} Let $(b_{ij})$ be a braiding of Cartan type  
and let us assume that $b_{ij}$ is a root of 1 of odd order for any $i, j$.  
Then there exists a finite abelian group $\Gamma$, sequences $g(1)$, \dots, 
$g(\theta) \in\Gamma$, $\chi(1)$, \dots, $\chi(\theta) \in\widehat{\Gamma}$ 
and $F \in \k \Gamma \otimes \k\Gamma$ such that 

\roster \item"(i)" for $1 \le i,j \le \theta$, we have 

$$b_{ij} = \langle \chi(j), g(i)\rangle; \tag 4.2$$ 

\item"(ii)" the braiding (1.7) corresponding to the Nichols algebra 
$R_{F}$ is symmetric and of Cartan type, with the same associated Cartan matrix and  diagonal
elements. \endroster 

Here $R$ is the Nichols algebra associated to the Yetter-Drinfeld module $V$ with  
basis $x_{1}, \dots, x_{\theta}$ and structure defined by (1.5).  \endproclaim

 \medpagebreak 
\demo{Proof} (i): We define $\Gamma := <y(1)> \oplus \dots \oplus <y(\theta)>$, 
where we impose that the order $E_{j}$ of $y_{j}$ is the least common multiple 
of the orders of $b_{\ell i}$, for all $i$, and for $\ell \leq j$. Then 
$E_{j}\vert E_{j+1}$: (2.1) holds. Let $q\in k$ be a root of 1 of order 
$E_{\theta}$. Let $D_{i} = \dsize\frac{E_{\theta}}{E_{i}}$, $1\le i \le \theta$ 
and let $q_{i} = q^{D_{i}}$. Then $q_{i}$ is a root of 1 of order $E_{i}$. 
 Let $\gamma(i)\in \widehat{\Gamma}$ be the unique character such that 
$\langle \gamma(i),  y(j)\rangle = q_{i}^{\delta_{ij}}$, for $1\le i,j\le 
N$. Then $\ord \gamma(\ell) = E_{\ell}$ and   $\widehat{\Gamma} = 
<\gamma(1)> \oplus \dots \oplus <\gamma(M)>$. We choose $\chi(j) = 
\gamma(j)$ for all $j$ and define $g(i)$ by (4.2); it makes sense 
because the order of $\langle \chi(j), g(i)\rangle$ divides $E_{j}$ 
for all $i,j$.    

\medpagebreak (ii): Let $\alpha_{ij}$ be given by (2.9). 
We are looking for integers $c_{ij}$ such that, if $F$ 
is associated to the cocycle $\omega$ given by (2.5), 
then the matrix $\alpha_{ij}^{F}$ is symmetric. Assume that $i < j$. 
By (2.10) and our preceding choices, we have 

$$\align \alpha_{ij}& =  D_{i}g(j)_{i}, \\ \alpha_{ij}^{F}&= -D_{i}c_{ij} + \alpha_{ij} 
=  -D_{i}c_{ij} + D_{i}g(j)_{i},\\ \alpha_{ji}^{F}&= D_{i}c_{ij} + \alpha_{ji} 
=  D_{i}c_{ij} + D_{j}g(i)_{j}; \endalign$$ hence $\alpha_{ij}^{F}
= \alpha_{ji}^{F}$ if and only if  $$2D_{i}c_{ij} = \alpha_{ij} - \alpha_{ji}
 = D_{i}g(j)_{i} - D_{j}g(i)_{j} \mod E_{\theta}. \tag 4.3$$ 

We solve (4.3). The factor 2 does not trouble because $E_{\theta}$ is odd. 
That is, it is enough to solve  

$$D_{i}\widetilde c_{ij}  = D_{i}g(j)_{i} - D_{j}g(i)_{j} \mod E_{\theta}. \tag 4.4$$ 

Now, since $g(j) = \sum_{t} g(j)_{t} y(t)$ in $\Gamma$, we have 

$$b_{ij} = \langle \chi(i), g(j)\rangle  = q_{i}^{g(j)_{i}} = q^{D_{i}g(j)_{i}}.$$ 

 As the braiding $(b_{ij})$ is of Cartan type, we see that $$q^{D_{i}g(j)_{i} + D_{j}g(i)_{j}} 
= b_{ij}b_{ji} = b_{ii}^{a_{ij}} = q^{D_{i}a_{ij}}$$ 
 and therefore 

$$D_{j}g(i)_{j} = D_{i}a_{ij} - D_{i}g(j)_{i}\mod E_{\theta}.$$ 

Thus $\widetilde c_{ij}  = 2g(j)_{i} - a_{ij}$ is a solution of (4.4).
The last statement follows from Remark 2.4 \qed \enddemo

In the rest of this section, we shall consider matrices ${\bold b} 
= (b_{ij})_{1\le i, j\le \theta}$ of Cartan type, such that  $b_{ij}$
 is a root of 1 of odd order and $b_{ij} = b_{ij}$ for all $i, j$. 
By Lemma 4.1, we can consider its associated Nichols algebra $R$; the 
comultiplication is an algebra map where the product in $R\otimes R$
 is determined by $$(x_{r}\otimes x_{i})(x_{j}\otimes x_{s}) = b_{ij} 
x_{r}x_{j}\otimes x_{i}x_{s}.$$ 

If $I\subset \{1, \dots, \theta\}$ then we denote by $R(I)$
 the Nichols algebra corresponding to the braiding $(b_{ij})_{i, j\in I}$. 
Clearly, there is an injective map of braided Hopf algebras $ R(I) \to R$.    

\proclaim{Lemma 4.2} (i). If $b_{ij} = 1$ then $x_{i}x_{j} = x_{j}x_{i}$.  

(ii). Assume there is $I\subset \{1, \dots, \theta\}$ such that $b_{ij} = 1$
 for  all $i\in I$, $j\in J := \{1, \dots, \theta\} - I$. 
Then  $R\simeq  R(I)\otimes  R(J)$. \endproclaim  

\demo{Proof} By a direct computation, we see that $x_{i}x_{j} - x_{j}x_{i}$
 is primitive; but the primitive elements are concentrated in degree one, 
so it should be 0. This proves (i).   Let $S =  R(I)\otimes  R(J)$; 
it is a graded algebra and coalgebra with respect to the grading 
$S(h) = \oplus_{\ell}  R(I)(\ell)\otimes  R(J)(h - \ell)$ and the usual 
tensor product multiplication and comultiplication. Then $S$ is a 
strictly graded coalgebra, {\it i.e.} (1.1), (1.2), (1.3) hold, 
{\it cf.} \cite{Sw, p. 240}. On the other hand, $S$ is a Yetter-Drinfeld 
module and it is easy to conclude from the hypothesis and (i) that it is 
a braided Hopf algebra.  Clearly, $S$ is generated as algebra by $S(1) 
\simeq V$, {\it i.e.} it satisfies (1.4). Thus $S$ is a Nichols algebra, hence 
isomorphic to $R$ because their spaces of primitive elements coincide. \qed\enddemo   

\bigpagebreak We consider the following equivalence relation on 
$\{1, \dots, \theta\}$: $i\sim j$ if there exists a sequence of 
elements $i = h_{0}, \dots, h_{P} = j$ in $\{1, \dots, \theta\}$ 
such that $h_{\ell} \neq h_{\ell + 1}$ and $b_{h_{\ell}  h_{\ell + 1}} 
\neq 1$. Note that $P = 0$ is allowed and gives $i\sim i$. A class of 
this relation is called a connected component. 
This is equivalent to the notion of connected Dynkin diagram or 
indecomposable matrix as in \cite{K, 1.7.1}.
The preceding Lemma 4.2 
shows that the Nichols algebra $R$ is the tensor product of the Nichols algebras corresponding 
to its connected components. Therefore, to prove Theorem 1.1 it is enough to 
 assume that 
$\bold b$ is {\it connected}.

\subsubhead 4.3 Conditions for FL-type\endsubsubhead

\proclaim{Lemma 4.3} Let ${\bold b} = (b_{ij})_{1\le i, j\le \theta}$ 
be a braiding of Cartan type such that  
$b_{ij}$ is a root of 1 of odd order and $b_{ij} = b_{ji}$ for all $i, j$. 

If $\bold b$ is symmetrizable and satisfies the relative primeness condition
($a_{ij}$ is 0 or relatively prime to the order of $b_{ii}$ for all $i,j$) 
then it is of FL-type. \endproclaim 

\demo{Proof} We can easily reduce to the case when $\bold b$ is connected.

We claim first that the order of $b_{kk}$ is the same for all
$k$. Indeed, let us fix $i\neq j$ such that $a_{ij} \neq 0$ (if no such pair $i,j$
exists, then by connectedness $\theta = 1$ and there is nothing to prove).
Then, since the orders of the entries are odd,
$$\ord b_{ij} = \ord b_{ij}^2 = \ord b_{ii}^{a_{ij}} = \ord b_{ii};$$
so $\ord b_{ii} = \ord b_{jj}$. The claim follows by connectedness.
We call $N = \ord b_{ii}$.

Let $d_{k}$ be integers as in (1.11). We can assume that they are relatively prime,
that is, $(d_{1}, \dots, d_{\theta}) = 1$. Then it is easy to see that
$N$ and $d_k$  are relatively prime for all $k$. 
Indeed, if $t$ divides $N$ and $d_i$ then $t$ divides $d_{j}$, 
since $t$  divides $d_ia_{ij} = d_{j}a_{ji}$ and $(N, a_{ji}) = 1$. Again, 
the claim follows by connectedness.

In particular, there is a unique root of unity $q$ of order $N$ such
that $b_{ii} = q^{2d_{i}}$. We claim finally that this $q$ satisfies (1.12). Indeed,
$b_{ii}^{a_{ij}} = q^{2d_{i}a_{ij}} = q^{2d_{j}a_{ji}}$ and $b_{ii}^{a_{ij}}
= b_{ij}^2 = b_{jj}^{a_{ji}}$. Thus, $ b_{ij} = q^{d_{i}a_{ij}}$ (since $N$
is odd) and  $b_{jj} = q^{2d_{j}}$ (since $(N, 2a_{ji}) = 1$). Once more, 
the claim follows by connectedness.
\qed\enddemo

\bigpagebreak
We now investigate conditions  for FL-type in rank 2. 
Let $\bold b = (b_{ij})_{1\le i,j \le 2}$ be a braiding of connected
Cartan type whose entries have odd order, such that $b_{12} = b_{21}$. 
Then it is automatically symmetrizable since the rank is 2.
"Connected" means that $a_{12} \neq 0$. Let $d_{1}$, $d_{2}$ be relatively
prime integers such that $d_{1} a_{12} = d_{2}a_{21}$. We denote by
$N_{i}$ the order of $b_{ii}$.

\proclaim{Lemma 4.4} The following are equivalent:

\roster\item"(i)" $\bold b$ is of FL-type.

\item"(ii)" There exists $u \in \k$ of odd order
such that $u^{2d_i} = b_{ii}$, $i = 1,2$.

\item"(iii)" There exists $v \in \k$ such that $v^{d_i} = b_{ii}$, $i = 1,2$.
\endroster \endproclaim 

\demo{Proof} (iii) $\implies$ (ii):  Note that the order of $v$ divides
$d_1N_1$ and $d_2N_2$. Since the $N_i$'s are odd and  $d_{1}$, $d_{2}$ 
are relatively prime, the order of $v$ is odd. Hence $v$ has a square root 
$u$ of odd order.

(ii) $\implies$ (i): We have 
$b_{12}^2 = b_{11}^{a_{12}} = u^{2d_{1}a_{12}}$. Hence $b_{12} = u^{d_{1}a_{12}}$ since
both have odd order. 

(i) $\implies$ (iii): If $q$ satisfies (1.12), take $v = q^2$.
\qed\enddemo

We want to give a criterion for the condition (iii) in the preceding Lemma.
 Let $e_{i}$ , $i = 1,2$ be non-zero
integers such that 
$$ e_{1}d_{1}N_{1} = e_{2}d_{2}N_{2} =: r;$$
for instance we could take $r$ the lowest common multiple of $d_{1}N_{1}$
and $d_{2}N_{2}$. Observe that there exists $s \in \Bbb Z$ such that $r= d_{1}d_{2}s$.

Let now $\xi \in \k$ be a primitive $r$-th root of 1 and choose integers
$k_1$, $k_2$ such that
$$ b_{ii} = \xi^{e_{i}d_{i}k_{i}};
$$
this is possible because $\xi^{e_{i}d_{i}}$ has order $N_i$. 

\proclaim{Lemma 4.5} Condition (iii) in Lemma 4.4  is equivalent to
\roster
\item"(4.5)" $e_{1}k_1 \equiv e_2k_2 \mod s.$
\endroster \endproclaim 

\demo{Proof} We claim first that (iii) is equivalent to the following statement:

\roster
\item"(iv)" There exists $t_i \in \Bbb Z$, $i=1,2$ such that 
$$\align e_{i}d_ik_i  &\equiv e_id_it_i \mod r, \quad i=1,2\tag 4.6 \\
e_{1}t_1 &\equiv e_2t_2 \mod r.\tag 4.7\endalign$$
\endroster
Indeed, if (iii) holds then $q^{d_iN_i} = 1$; as $\xi^{e_{i}}$ 
has order $d_{i}N_i$, there exists $t_i \in \Bbb Z$ such that
$q = \xi^{e_{i}t_i}$. Then (4.6), (4.7) follow now without difficulty.
Conversely, if (iv) holds, take $q = \xi^{e_{1}t_{1}} = \xi^{e_{2}t_{2}}$,
by (4.7). Then (iii) is true by (4.6). The claim is proved.

Now (4.6) is equivalent to $ k_i  \equiv t_i \mod N_i$, $i=1,2$.

Assume (iv). Then there exist $x_i \in \Bbb Z$ such that 
$t_{i} = k_{i} + N_{i}x_{i}$, $i=1,2$. Now $r$ (and {\it a fortiori} $s$)  divides
$$
e_{1}t_1 - e_2t_2 = e_{1}k_1 - e_2k_2 +e_{1}N_{1}x_{1} - e_{2}N_{2}x_{2};$$
but $e_{1}N_{1} = d_{2}s$ and $e_{2}N_{2} = d_{1}s$. So (4.5) holds.

If (4.5) holds, let $y\in \Bbb Z$ such that $e_{1}k_1 - e_2k_2 = -ys$.
As $d_{1}$ and $d_{2}$ are relatively prime, there exists $x_{1}, x_{2}$
such that $y = d_2x_1 - d_1x_2$. If now we take $t_{i} = k_{i} + 
N_{i}x_{i}$, $i=1,2$, then (4.6) holds by definition; and
$$
e_{1}t_1 - e_2t_2 = e_{1}k_1 - e_2k_2 +e_{1}N_{1}x_{1} - e_{2}N_{2}x_{2}
= e_{1}k_1 - e_2k_2 +  d_{2}sx_{1} - d_{1}sx_{2} = 0.$$
So, (4.7) follows. \qed\enddemo

\medpagebreak
Lemma 4.5 describes an easy algorithm to decide whether a given $\bold b$
is locally of FL-type. As an example we note:

\proclaim{Corollary 4.6}  There are braidings of symmetrizable Cartan type
which are not of FL-type. 
\endproclaim 
\demo{Proof} Let $p$ be an odd prime number and
take $N_{1} =  N_{2} = p^2$,  $a_{12}= a_{12}= - p$. Let $k_{1}$, $k_{2}$
be two elements not divisible by $p$ such that
$$ k_{1}\not\equiv k_{2} \mod p^2, \qquad  k_{1} \equiv k_{2} \mod p.
$$
Let $q$ be a root of unity of order $p^2$ and let $b_{ii} = q^{k_i}$, 
$i=1,2$, $b_{12} = b_{21}$ the unique root of unity of odd order such that 
$b_{12}^2= b_{11}^p$. Then  $b_{12}^2= b_{22}^p$, but $s = p^2$
and $e_{1} = e_{2} = 1$, so (4.5) does not hold. \qed\enddemo

\medpagebreak
\proclaim{Example 4.7}There are braidings of  Cartan type
which are locally of FL-type but not symmetrizable. 
\endproclaim 

\demo{Proof} Take $\bold b = (q^{d_{i}a_{ij}})$, where $q$ is a root of 1
of order 5, $d_{1} = 1$,  $d_{2} = 2$, $d_{3} = 3$ and 
$$
(a_{ij}) = \pmatrix 2 & -3 & -3 \\  -4 & 2 & -3 \\ -1 & -2 & 2\endpmatrix. 
$$ 
Then $\bold b$ is of Cartan type since $d_{i}a_{ij} = d_{j}a_{ji} \mod 5$
for all $i,j$. By Lemma 4.3,  $\bold b$ is locally of FL-type.
But $(a_{ij})$ is not symmetrizable. \qed\enddemo

\medpagebreak
\proclaim{Lemma 4.8} Let $\bold b$ be a braiding of Cartan type and rank 2. 
If $a_{12} = -1$ and $a_{21}$ is odd, 
then $\bold b$ is of FL-type. 
\endproclaim 
\demo{Proof} Let $n = -a_{21}$. Let $\widetilde q$ be a root of 1 such that 
$\widetilde q^{n} = b_{12}^{-1}$; then  $b_{22}^{n} = b_{11} = 
b_{12}^{-2} = \widetilde q^{2n}$.  Hence there exists an $n$-th root of 
1 $\omega$ such that $b_{22} = \omega^{2}\widetilde q^{2}$ since $n$ is odd.  
We choose  then $q = \omega\widetilde q$. 
\qed\enddemo

\bigpagebreak
\demo{Proof of Theorem 1.1} 

\medpagebreak

By Lemma 4.1 we can assume that  $\bold b$ is symmetric.
By Lemma 4.2 we can assume that  $\bold b$ is connected.

\medpagebreak
Part (1). By Lemma 4.3 (if the type is different from $G_2$)
or Lemma 4.8 (if the type is $G_2$), the braiding $\bold b$ 
is of FL-type. We now apply Theorem 3.1.

\medpagebreak
Part (2). For any fixed pair $i,j$ in 
$\{1, \dots, \theta\}$, the  braiding of rank 2 corresponding to the submatrix
supported by $i,j $ is of FL-type and $\goth B(V)$
is finite dimensional. By Theorem 3.1, we conclude that $a_{ij}a_{ji} = 0, 1, 2$ or $3$.
Since 3 is relatively prime to the orders of the $b_{ii}$ whenever $a_{ij} = -3$ 
for some $j$ by hypothesis,
we can apply Lemma 4.3 and conclude again from Theorem 3.1 that $\bold b$
is of finite type whenever it is symmetrizable.

\medpagebreak
So, it only remains to show that: $\bold b$ is necessarily symmetrizable.
   
 It is known that a matrix
which is either  simply-laced or has  no cycles is symmetrizable, cf. \cite{K, ex. 2.1}. 
If the matrix has a cycle, its corresponding Nichols algebra   is still
finite dimensional.   
We are reduced to prove that 
no cycle with a double or triple arrow can arise.  

Assume first that $\bold b$ is a cycle with a triple arrow. If $\theta \ge 4$,
then we remove a suitable vertex and get a subdiagram with no cycles
which is not of finite type. This is a contradiction.

So, assume $\theta = 3$. Then there are several possibilities of cycles with
(at least) one triple arrow. We can discard two of them which are symmetrizable.
We discard the rest because of the restrictions on the orders of $b_{ii}$. 
Let us proceed with the details. Consider
a Cartan matrix

$$
(a_{ij}) = \pmatrix 2 & -3 & -b \\  -1 & 2 & -c \\ -d & -e & 2\endpmatrix.
$$ 
Here $b,c,d,e$ are positive integers such that $bd, ce = 1,2$ or 3.
This means that 
$$
b_{11}^3 = b_{22}, \quad b_{22}^c = b_{33}^e, \quad b_{11}^b = b_{33}^d.
$$
Hence, 
$$
 b_{11}^b = b_{33}^d, \quad  b_{11}^{3c} = b_{33}^e.
$$ 
We know that $b = 1$ or $d=1$. Assume $b=1$. Then $b_{33}^{3cd - e} = 1$. 
We have several subcases according to the values of $c,d,e$.

\roster
\item"(i)" If $e = 3$, then $c = 1$ and  $b_{33}^{3(d - 1)} = 1$. 
There are three possibilities:
\endroster

If $d=1$ the matrix is symmetrizable and we can apply Theorem 3.1.
If $d=2$,  $b_{33}^{3} = 1$. 
If $d=3$,  $b_{33}^{6} = 1$. As we assume that the orders are odd, we should
have $b_{33}^{3} = 1$. These two possibilites are not possible by our choice 
following the rule (1.10).

\roster
\item"(ii)" If $e = 2$, then $c = 1$ and  $b_{33}^{3d - 2} = 1$. 
There are three possibilities:
\endroster

If $d=1$ $b_{33} = 1$, a possibility that we excluded.
If $d=2$,  $b_{33}^{4} = 1$, but we assume that the orders are odd. 
If $d=3$,  $b_{33}^{7} = 1$. For this reason we exclude the order 7.

\roster
\item"(iii)" If $e = 1$, then    $b_{33}^{3cd - 1} = 1$. 
There are various possibilities:
\endroster

As $cd = 1,2,3,4,6,9$, $3cd - 1$ could take the values $2,5,8,11,17,26$. 
Since we assume that the orders are odd,  we have to exclude the orders 5, 11, 13, 17.

Let us now assume that $d = 1$. Then $b_{11}^{3c - be} = 1$. 

\roster
\item"(iv)" If $c = 3$, then $e = 1$ and  $b_{11}^{9 - b} = 1$. 
There are three possibilities:
\endroster

If $b=1$,  $b_{11}^{8} = 1$; but the orders are odd.
If $b=2$,  $b_{11}^{7} = 1$. This is a new instance where 7 should be excluded. 
If $b=3$,  $b_{11}^{6} = 1$. As we assume that the orders are odd, we should
have $b_{11}^{3} = 1$. This contradicts the rule (1.10).

\roster
\item"(v)" If $c = 2$, then $e = 1$ and  $b_{11}^{6 - b} = 1$. 
There are three possibilities:
\endroster

If $b=1$,  $b_{11} = 5$. This is a new instance where 5 should be excluded.
If $b=2$,  $b_{11}^{4} = 1$, but we assume that the orders are odd. 
If $b=3$,  $b_{11}^{3} = 1$.  This contradicts the rule (1.10).

\roster
\item"(vi)" If $c = 1$, then    $b_{11}^{3 - be} = 1$. 
There are various possibilities:
\endroster

As $be = 1,2,3,4,6,9$, $3 - be$ could take the values $2, 1, 0, -1, -3, -6$. 
The values $2,1,-1$ are excluded by hypothesis; the values $-3, -6$ are excluded
by the rule (1.10). The value 0 arises only when $\bold b$ is symmetrizable,
but then we can apply Theorem 3.1.

We next assume that $\bold b$ is a cycle with no triple arrow but
at least a double one. If $\theta \ge 6$, again we have a contradiction by 
removing a suitable vertex.  If $\theta = 5$, we remove systematically
different vertices. Since the resulting subdiagrams are symmetrizable and hence
of finite type, we see that $\bold b$ has only one double arrow. We 
can easily get a contradiction.

Now we assume $\theta = 4$, with a Cartan matrix
$$
(a_{ij}) = \pmatrix 2 & -2 & 0 & -b \\  -1 & 2 & -c & 0 \\ 0 & -d & 2 & -e
\\ -f & 0 & -g  & 2\endpmatrix.
$$
Here $b,c,d,e, f, g$ are positive integers such that $bf, cd, eg = 1,2$. By
removing the third, respectively the fourth,
 vertex, we conclude that $bf = cd = 1$. There are 
three possibilities: $(e,g) = (1,1), (1,2)$ or $(2,1)$. 
The first is not possible by a similar argument as for $\theta = 5$. The second corresponds
to a symmetrizable matrix and the third is the first instance forcing to exclude the order 3.

We finally assume that $\theta = 3$. We have a 3-cycle with no triple arrow but
at least a double one. Arguing as in the case of triple arrows, we eliminate all the possibilities
except two; for those we exclude the orders 3 and 7. For instance,
$$
(a_{ij}) = \pmatrix 2 & -2 & -1 \\  -1 & 2 & -2 \\ -2 & -1 & 2\endpmatrix, 
$$
then 
$$
b_{11}^2 = b_{22}, \quad b_{22}^2 = b_{33}, \quad b_{11} = b_{33}^{2}.
$$
Hence $b_{11} = b_{11}^8$, which is not possible because  the order of $b_{11}$ is not 7.
\qed\enddemo

\remark{Remark}
We observe that the preceding proof of Theorem 1.1 Part (2) only needs
the hypothesis "locally of FL-type" to conclude that 
\roster\item"(4.8)" A finite dimensional
Nichols algebra 
$\goth B(V)$ of rank 2 is necessarily of finite type. \endroster

If $\Gamma$
is any finite abelian group such that (4.8) holds for $\goth B(V)$
in $^{\Gamma}_{\Gamma}{\Cal YD}$ then the conclusion of Theorem 1.1 Part (2)
is valid for all  $\goth B(V)$
in $^{\Gamma}_{\Gamma}{\Cal YD}$.
\endremark

\bigpagebreak
For further use, we list the Cartan matrices, up to numbering,
 causing troubles for
small values of $p$:

\roster
\item"(4.9)" For $p=3$, the matrix
$$
(a_{ij}) = \pmatrix 2 & -2 & -1 \\  -1 & 2 & -1 \\ -2 & -1 & 2\endpmatrix.
$$

\item"(4.10)" For $p=3$, the matrix
$$
(a_{ij}) = \pmatrix 2 & -2 & 0 & -1 \\  -1 & 2 & -1 & 0 \\ 0 & -1 & 2 & -2
\\ -1 & 0 & -1  & 2\endpmatrix.
$$

\item"(4.11)" For $p=5$, the matrix
$$
(a_{ij}) = \pmatrix 2 & -3 & -1 \\  -1 & 2 & -1 \\ -2 & -1 & 2\endpmatrix.
$$ 

\item"(4.12)" For $p=5$, the matrix
$$
(a_{ij}) = \pmatrix 2 & -3 & -1 \\  -1 & 2 & -2 \\ -1 & -1 & 2\endpmatrix.
$$

\item"(4.13)" For $p=7$, the matrix
$$
(a_{ij}) = \pmatrix 2 & -3 & -1 \\  -1 & 2 & -1 \\ -3 & -2 & 2\endpmatrix.
$$ 

\item"(4.14)" For $p=7$, the matrix
$$
(a_{ij}) = \pmatrix 2 & -2 & -1 \\  -1 & 2 & -2 \\ -2 & -1 & 2\endpmatrix.
$$

\item"(4.15)" For $p=11$, the matrix
$$
(a_{ij}) = \pmatrix 2 & -3 & -1 \\  -1 & 2 & -2 \\ -2 & -1 & 2\endpmatrix.
$$

\item"(4.16)" For $p=13$, the matrix
$$
(a_{ij}) = \pmatrix 2 & -3 & -1 \\  -1 & 2 & -3 \\ -3 & -1 & 2\endpmatrix.
$$ 

\item"(4.17)" For $p= 17$, the matrix
$$
(a_{ij}) = \pmatrix 2 & -3 & -1 \\  -1 & 2 & -2 \\ -3 & -1 & 2\endpmatrix.
$$

\endroster

 \bigpagebreak \subhead \S 5. Nichols algebras over $\Bbb Z/(p)$\endsubhead

In this section, $\Gamma$ will denote $\Bbb Z/(p)$, $p$ an odd prime.
 We first discuss Nichols algebras of rank 2. If $V$ is a Yetter-Drinfeld module
 of dimension 2 satisfying (1.6), then there exists a generator $u$ of 
$\Gamma$, $q \in k^{\times}$ of order $p$ and integers $b$, $d$ such that 

$$ g(1) = u, \quad g(2) = u^{b}, \quad \langle \chi(1),  u\rangle 
= q, \quad  \langle \chi(2),  u\rangle = q^{d}. \tag 5.1$$  

So $\dsize b_{11} = q, \quad b_{22} = q^{bd}$.

\example{Dynkin diagram of type $A_{2}$}  We should have 

$$\langle\chi(1),  g(2)\rangle \langle \chi(2), 
 g(1)\rangle = b_{11}^{-1} = b_{22}^{-1}. $$ 

This means $b + d \equiv -1 \equiv - bd \mod p$. 
It follows readily that $b^{2} + b + 1\equiv 0 \mod p$. 
Hence $p = 3$ and $b \equiv 1 \mod p$, or $b \not \equiv 1\mod p$ 
is a cubic root of 1. In the last case, $b$ exists if and only if 
$p-1$ is divisible by 3. The corresponding matrices are not symmetric.   \endexample  

\bigpagebreak \example{Dynkin diagram of type $B_{2}$} We should have 

$$ \langle\chi(1),  g(2)\rangle \langle \chi(2),  g(1)\rangle = 
b_{11}^{-1} = b_{22}^{-2}. $$ 

This means $b + d \equiv -1 \equiv - 2bd \mod p$. 
Thus $2b^{2} + 2b + 1 \equiv 0 \mod p$ and looking at the discriminant of 
this equation we see that it has a solution if and only if $-1$ is a 
square, {\it i.e.} exactly when $p \equiv 1 \mod 4$. The corresponding 
matrices are not symmetric.\endexample  

\bigpagebreak \example{Dynkin diagram of type $G_{2}$} By (1.10), we have $p > 3$. 
We should have 
$$ \langle\chi(1),  g(2)\rangle \langle \chi(2),  g(1)\rangle = b_{11}^{-1} 
=  b_{22}^{-3}. $$ 

This means $b + d \equiv -1 
\equiv - 3bd \mod p$. Thus $3b^{2} + 3b + 1 \equiv 0 \mod p$ and looking 
at the discriminant of this equation we see that it has a solution if and 
only if $-3$ is a square $\mod p$. By the quadratic reciprocity law, 
this happens exactly when $p \equiv 1 \mod 3$.   The matrices are not 
symmetric. \endexample  

\bigpagebreak We now conclude the proof of the existence part of Theorem 1.3.

\proclaim{Proposition 5.1} There is no Nichols algebra of rank  $\geq 3$ over 
$\Bbb Z/(p)$ with finite Cartan matrix except when $p = 3$ and  the 
corresponding Dynkin diagram is $A_{2}\times A_{1}$ or $A_{2}\times 
A_{2}$.\endproclaim 

\demo{Proof}We first consider Dynkin diagrams of rank 3. We have the following cases: 

$$\alignat3 (a)&\, A_{2}\times A_{1}, &\qquad (b)&\, B_{2}\times A_{1}, 
&\qquad (c)&\, G_{2}\times A_{1},\\ (d)&\, A_{3}, &\qquad (e)&\, B_{3}, 
&\qquad (f)&\, C_{3}, \qquad (g) \, A_{1}\times A_{1}\times A_{1}.
\endalignat $$ 

The case (g) does not arise, as shown in \cite{AS2, Section 4}. 
Therefore we can fix a numeration of the vertices of the Dynkin diagram 
such that the vertices 1 and 2 (resp., 1 and 3) are connected (resp., 
not connected). We also assume that there is only one arrow between 1 and 2 
in cases (e) an (f).  That is, 

$$\alignat3 (a)&\, \underset{1}\to\circ \frac{\quad}{\quad} 
\underset{2}\to\circ \qquad \underset{3}\to\circ, &\qquad  
(b)&\,\underset{1}\to\circ \Rightarrow \underset{2}\to\circ \qquad 
\underset{3}\to\circ, &\qquad  (c)&\, \underset{1}\to\circ \Rrightarrow 
\underset{2}\to\circ \qquad \underset{3}\to\circ,\\  (d)&\, \underset{1}
\to\circ \frac{\quad}{\quad} \underset{2}\to\circ \frac{\quad}{\quad} 
\underset{3}\to\circ, &\qquad  (e)&\, \underset{1}\to\circ \frac{\quad}
{\quad} \underset{2}\to\circ \Rightarrow \underset{3}\to\circ, &\qquad  
(f)&\, \underset{1}\to\circ \frac{\quad}{\quad} \underset{2}\to\circ 
\Leftarrow \underset{3}\to\circ.\endalignat $$    

If $V$ is a Yetter-Drinfeld module of dimension 3 satisfying (1.6), 
then there exist a generator $u$ of $\Gamma$, $q \in k^{\times}$ of 
order $p$ and integers $b, d, e, f$ (none of them divisible by $p$) 
such that   

$$\alignat3 g(1) & = u, & \quad g(2) & = u^{b}, & \quad g(3) & = u^{e}, 
\\ \langle \chi(1),  u\rangle & = q, & \qquad  \langle \chi(2),  u\rangle 
& = q^{d},  & \qquad \langle \chi(3),  u\rangle&  = q^{f}. \endalignat$$ 

 So $\dsize b_{11} = q$, $b_{22} = q^{bd}$, $b_{33} = q^{ef}$. Also,  
$ \langle\chi(1),  g(3)\rangle \langle \chi(3),  g(1)\rangle = 1$ means  

$$e + f \equiv 0 \mod p. \tag 5.2$$ 

Now, considering the subdiagram supported by 
the vertices 1 and 2, we conclude from the arguments above for the 
rank 2 case that   

\roster   \item"(5.3)" In cases (a), (d), (e) and (f), 
$b^{2} + b + 1\equiv 0 \mod p$, $bd \equiv 1 \mod p$ and $p= 3$ or 
$p\equiv 1 \mod 3$.   

\item"(5.4)" In case (b), 
$2b^{2} + 2b + 1\equiv 0 \mod p$, $2bd \equiv 1 \mod p$ and $p\equiv 1 
\mod 4$.  

\item"(5.5)" In case (c), $3b^{2} + 3b + 1\equiv 0 \mod p$, 
$3bd \equiv 1 \mod p$ and $p\equiv 1 \mod 3$. \endroster  

On the other hand, in cases (a), (b) and (c) we have $ \langle\chi(2),  
g(3)\rangle \langle \chi(3),  g(2)\rangle = 1$, {\it i.e.}  

$$ed + bf \equiv 0 \mod p. \tag 5.6$$ 

Now, combining (5.6) with (5.2) we conclude that $b \equiv d \mod p$.  Hence:   

\medpagebreak In case (a), $b^{2} \equiv 1 \mod p$. But also $b^{3} \equiv 1  \mod p$, 
so $b \equiv 1 \mod p$ and $p = 3$.  

\medpagebreak In case (b), $2b^{2} \equiv 1 \mod p$ and thus $b \equiv -1 \mod p$ and 
$1 \equiv 0 \mod p$, a contradiction.  

\medpagebreak In case (c), $3b^{2} \equiv 1 \mod p$. Plugging this into the first 
equation of (5.5) we easily get a contradiction.

\medpagebreak
Now we turn to the remaining cases. In case (d),  we have 
$$ \langle\chi(2),  g(3)\rangle \langle \chi(3),  g(2)\rangle = 
b_{22}^{-1} =b_{33}^{-1};$$ using (5.2), this implies 

$$f(b - d) \equiv -1 \equiv f^{2} \mod p. $$ 

Thus $b - d \equiv f\mod p$ and $(b - d)^{2} \equiv b^{2} + b - 2 
\equiv - 1\mod p$; this last equation contradicts the first of (5.3).  
In case (e),  we have $\langle\chi(2),  g(3)\rangle \langle \chi(3),  
g(2)\rangle = b_{22}^{-1} =b_{33}^{-2}$; using (5.2), this implies 

$$f(b - d) \equiv -1 \equiv 2f^{2} \mod p. $$ 

Thus $b - d \equiv 2f\mod p$ and $(b - d)^{2} \equiv b^{2} + b - 2 
\equiv  4f^{2} \equiv - 2\mod p$; this last equation contradicts the first 
of (5.3).  In case (f),  we have $\langle\chi(2),  g(3)\rangle \langle 
\chi(3),  g(2)\rangle = b_{22}^{-2} =b_{33}^{-1}$; using (5.2), 
this implies 

$$f(b - d) \equiv -2 \equiv f^{2} \mod p. $$ 

Thus $b - d \equiv f\mod p$ and $(b - d)^{2} \equiv b^{2} + b - 2 
\equiv  f^{2} \equiv - 2\mod p$; this last equation contradicts the 
first of (5.3).   

\medpagebreak We have shown that the only possibility in rank 3 is case (a) 
with $p = 3$. Then, if $V$ is a Yetter-Drinfeld module of Cartan type satisfying 
(1.6) with finite Cartan matrix and rank $\geq 3$, then $p = 3$, the 
corresponding Dynkin diagram should at most 2 connected components-- 
{\it cf.} case (g)--   and each component is of type $A_{1}$ or $A_{2}$. 
On the other hand, let $u$ be a generator of $\Bbb Z/(3)$ and $q$ a root of 1
 of order 3. It is easy to see that the sequences  

$$\aligned  g(1) &= u = g(2), \quad g(3) = u^{e} = g(4), \\ \langle \chi(1),  u\rangle& = q 
=  \langle \chi(2),  u\rangle, \quad \langle \chi(3),  u\rangle = q^{-e} = \langle \chi(4),  
u\rangle \endaligned \tag 5.7$$ 

\noindent define a Yetter-Drinfeld module of rank 4 over $\Bbb Z/(3)$, which is of 
Cartan type and has Dynkin diagram $A_{2}\times A_{2}$. \qed\enddemo

\proclaim{Lemma 5.2} Let $\goth B(V)$ be a finite dimensional Nichols algebra  over 
$\Bbb Z/(p)$. Then $V$ is of  finite Cartan type.\endproclaim

\demo{Proof} We know that the Lemma is true whenever $p$ is different from
3,5,7,11,13,17. Assume that $p$ is one of these small primes.
For all of them, it would be enough to prove that no Cartan
matrix as in (4.9), \dots, (4.17)  is possible over $\Bbb Z/(p)$.
All these matrices have a subdiagram of rank 2 of type $B_2$ or $G_2$.
This eliminates all the primes except for 13. We show then that no matrix
like (4.16) exists over $\Bbb Z/(13)$. Suppose in the contrary that
it exists. Then there exists a generator $u$ of $\Gamma$, $q \in k^{\times}$ of 
order $p$ and integers $b, c, d, e$ (none of them divisible by $13$) 
such that  $ g(1)  = u$, $g(2)  = u^{b}$, $g(3)  = u^{c}$, 
$\langle \chi(1),  u\rangle = q$, $\langle \chi(2),  u\rangle 
 = q^{d}$,  $\langle \chi(3),  u\rangle  = q^{e}$; and the following equations hold

$$
\align
b + d \equiv -3 \equiv - bd \mod 13, \tag 5.8 \\
c + e \equiv -1 \equiv - 3ec \mod 13, \tag 5.9 \\
be + dc \equiv -ec \equiv - 3bd \mod 13. \tag 5.10 
\endalign
$$
Now (5.8) implies $(b,d) = (2,8)$ or $(8,2)$; and 
(5.9) implies $(c,e) = (5,7)$ or $(7,5)$. We have four possibilities for 
$be + dc$ but none of them gives $-9 \mod 13$; that is, (5.10) does not hold.
\qed\enddemo

\proclaim{Lemma 5.3}   Let $R$ be a finite dimensional Hopf algebra
braided Hopf algebra in $_{\Gamma}^{\Gamma}\Cal{YD}$, $\Gamma = \Bbb Z/(p)$.
Then $\dim P(R) \le 2$ for $p>3$; and $\dim P(R) \le 4$ for $p = 3$.
\endproclaim

\demo{Proof}Consider first the coradical filtration of $R$ and the subalgebra $S$ of
the corresponding graded coalgebra generated by $P(R)$. Then $S$ is a Nichols algebra
and the claim follows from Proposition 5.1 and  Lemma 5.2. \qed\enddemo

\proclaim{Proposition 5.4} Let $R$ be a finite dimensional 
braided Hopf algebra in $_{\Gamma}^{\Gamma}\Cal{YD}$, $\Gamma = \Bbb Z/(p)$.
Assume that (1.1), (1.2) and (1.3) hold. Then $R$ is a Nichols algebra.
\endproclaim
\demo{Proof} Let $S = R^*$; it is a braided Hopf algebra in $_{\Gamma}^{\Gamma}\Cal{YD}$
and satisfies (1.1), (1.2) and (1.4); that is, $S$ is graded, $S_0 = \k 1$ and is 
generated in degree 1. We want to show that (1.3) holds for $S$; hence
(1.4) holds for $R$. See Lemma 8.1 below. 

We know that $S(1) \subseteq P(S)$ and we want to prove the equality.
If $\dim S(1) = 1$, we are done by \cite{AS2, Th. 3.2}.
If $\dim S(1) \ge 2$ and $p >3$ then we are done by Lemma 5.3. We can then assume that $p=3$; this case is treated in Section 8. \qed\enddemo
 
 \bigpagebreak \subhead \S 6. Isomorphisms between bosonizations\endsubhead

We begin with a general Lemma whose proof is straightforward. We shall denote 
the comultiplication of $R$ by $\Delta_{R}(r) = \sum r^{(1)}\otimes r^{(2)};$ 
we omit most of the time the summation sign.  

\proclaim{Lemma 6.1} Let $H$ be a Hopf algebra, $\psi: H\to H$ 
an automorphism of Hopf algebras, $V, W$  Yetter-Drinfeld modules over 
$H$.  

 \roster \item Let $V^{\psi}$ be the same space underlying $V$ but 
with action and coaction 
$$h._{\psi}v = \psi(h).v, \quad \delta^{\psi} (v)= \left(\psi^{-1} 
\otimes \id\right)\delta(v), \quad h\in H, v\in V.  \tag 6.1$$  

Then $V^{\psi}$ is also a Yetter-Drinfeld module over $H$. If $T: V \to W$ 
is a morphism in $^{H}_{H}{\Cal YD}$, then $T: V^{\psi} \to W^{\psi}$ also 
is. Moreover, the braiding $c: V^{\psi}\otimes W^{\psi} \to W^{\psi}
\otimes V^{\psi}$ coincides with the braiding $c: V\otimes W \to W\otimes V$.   

\medpagebreak \item If $R$ is an algebra (resp., a coalgebra, a Hopf algebra) 
in $^{H}_{H}{\Cal YD}$ then $R^{\psi}$ also is, with the same structural
 maps.   

\medpagebreak \item Let $R$ be a Hopf algebra in $^{H}_{H}\Cal{YD}$. 
Then the map $\Psi: R^{\psi}\# H \to R \# H$ given by $\Psi(r\# h) = r\# 
\psi(h)$ is an isomorphism of Hopf algebras.\endroster \endproclaim 

\demo{Proof}The proof of (1) is a direct computation that we leave to the 
reader.  The statement (2) follows from (1). Let us check (3). 
If $r, s\in R$ and $h,g\in H$ then 

$$\align \Psi\left((r\#h) (s\#g)\right)& 
=  \Psi\left( r h_{(1)}._{\psi}s \# h_{(2)}g \right) = r \psi(h_{(1)}).s 
\# \psi(h_{(2)}g);\\  \Psi(r\#h) \Psi(s\#g)& =  (r\#\psi(h)) (s\#\psi(g)) 
= r \psi(h)_{(1)}.s \# \psi(h)_{(2)}\psi(g).\endalign $$ 

That is, $\Psi$ is an algebra map. On the other hand, 

$$\align \Delta\Psi(r\#h) & =  \Delta(r\# \psi(h))  
= r^{(1)}\#  (r^{(2)})_{(-1)} \psi(h)_{(1)} \otimes 
(r^{(2)})_{(0)} \# \psi(h)_{(2)};\\    
(\Psi\otimes \Psi)\Delta(r\#h)& =  (\Psi\otimes \Psi) 
\left(r^{(1)}\#  \psi^{-1}\left((r^{(2)})_{(-1)}\right)h_{(1)} 
\otimes (r^{(2)})_{(0)} \# h_{(2)}\right) \\ &= r^{(1)}
\# (r^{(2)})_{(-1)} \psi(h_{(1)}) \otimes (r^{(2)})_{(0)} 
\# \psi(h_{(2)}).\endalign $$ 

Thus $\Psi$ is a bialgebra map and {\it a fortiori} a Hopf algebra map. \qed\enddemo     

\bigpagebreak Let $H$ be a cosemisimple Hopf algebra. 
Let  $A = \oplus_{n\geq 0} A(n)$ be a coradically graded Hopf algebra 
with coradical $A_{0} = A(0)$ isomorphic to $H$. Let $\iota: A(0) \to A$ 
denote the inclusion. Let $\pi: A\to A(0)$ be the unique graded projection 
with image  $A(0)$ and let $R = \{a\in A: (\id\otimes \pi) \Delta(a) = 
a\otimes 1\}$. Let $B$ be another Hopf algebra with coradical $B_{0}$ 
isomorphic to $H$.   

\proclaim{Lemma 6.2}  Let $\Phi: A\to B$ be an 
isomorphism of Hopf algebras.  Let $B(n) := \Psi\left(A(n) \right)$; 
then $B$ is also coradically graded with respect to the grading $B = \oplus_{n\geq 0} B(n)$. Let $\psi: H \to H$ denote the isomorphism induced
 by the restriction of $\Psi$, with respect to fixed identifications 
$A_{0}\simeq H$, $B_{0}\simeq H$.  Let $\jmath$, $\zeta$, $S$ 
denote the same objects as $\iota$, $\pi$, $R$ but with respect to the grading in $B$. 
Then $\Psi(R) = S$ and the restriction $\phi: R \to S$ of $\Psi$ 
is a morphism $R\to S^{\psi}$  of braided Hopf algebras in 
$_{H}^{H}\Cal{YD}$. Moreover, the following diagram commutes 

$$ \CD A@>\Psi>> B \\ @V\simeq VV @V\simeq VV \\ R\# H @>\phi \otimes  
\psi >> S\# H \endCD$$ 

where the vertical isomorphisms are the inverses of the bosonization maps.  
\endproclaim   

Note that $R\# H @>\phi \otimes  \psi >> S\# H$ factorizes as the composition
 of two Hopf algebra maps: $R\# H @>\phi \otimes\id>> S^{\psi}\# H @>
\id\otimes  \psi >> S\# H$.  \demo{Proof}  It is clear that 
$B = \oplus_{n\geq 0} B(n)$ is a grading of Hopf algebras. Since 
$\Psi$ preserves the terms of the coradical filtration, we have 

$$ B_{m} = \Psi(A_{m}) = \Psi\left(\oplus_{0\leq n \leq m}A(n)\right) 
= \oplus_{0\leq n \leq m}B(n); $$ that is, $B$ is coradically graded. 
It follows at once that 
$$\jmath\zeta = \Psi\iota \pi\Psi^{-1}.\tag 6.2$$ 

Hence $\Psi(R) = S$. We check that $\phi$ is a morphism from $R$ to 
$S^{\psi}$. It preserves the multiplication  because is the restriction 
of the algebra map $\Psi$. It preserves the comultiplication  because  
$\Psi$ is a coalgebra map and by (6.2). 
It preserves the action and the coaction: 

$$\align \phi(h.r) &= \Psi\left(h_{(1)} r \Ss(h_{(2)}) \right) =  
\psi\left(h_{(1)}\right) \phi(r) \psi\left(\Ss(h_{(2)})\right) 
\\ &= (\psi(h))_{(1)} \phi(r) \Ss\left((\psi(h))_{(2)}\right) = 
\psi(h).\phi(r) = h._{\psi} \phi(r); \\ \delta\phi(r) 
&= (\zeta\otimes \id) \Delta \Psi(r) = (\zeta\otimes \id) 
(\Psi\otimes \Psi) \Delta (r)  \\ &=  (\id\otimes \phi) 
(\psi^{-1} \pi\otimes \id) \Delta(r) = (\id\otimes \phi) \delta(r). 
\endalign  $$ 

Here again we used (6.2). The last statement follows because $\Psi$ is 
an algebra map. \qed\enddemo     

Let now $H = k\Gamma$, where $\Gamma$ is a finite abelian group. 
Let $g(1)$, \dots, $g(\theta) \in\Gamma$, $\chi(1)$, \dots, 
$\chi(\theta) \in\VGamma$ and let $V$ be the Yetter-Drinfeld module 
with structure given by (1.5). Let also $h(1)$, \dots, $h(\theta) 
\in\Gamma$, $\eta(1)$, \dots, $\eta(\theta) \in\VGamma$ and let 
analogously $W$ be the Yetter-Drinfeld module with structure given by 
(1.5) with $g(i)$'s, resp. $\chi(j)$'s, replaced by $h(i)$'s, resp. 
$\eta(j)$'s. Let $R$, $S$ be the corresponding Nichols algebras and $A = R\# H$, 
$B = S\# H$.   

\proclaim{Proposition 6.3} Assume that $R$ and $S$ are 
finite dimensional. The Hopf algebras $A$ and $B$ are isomorphic 
if and only if there exist $\varphi\in \Aut (\Gamma)$ and $\sigma \in \Bbb 
S_{\theta}$ such that 

$$g(j) = \varphi^{-1}\left(h(\sigma j)\right), 
\qquad  \chi(j) = \eta(\sigma j)\varphi, \qquad 1\leq j \leq \theta.  \tag 6.3$$
 \endproclaim  

\demo{Proof} Assume that $\Psi: A \to B$ is an isomorphism of Hopf algebras.  
We claim that $\Psi(A(n)) = B(n)$ for all $n \ge 0$. This is clear for 
$n=0$ since $\Psi$ preserves the coradical filtration. For $n = 1$ we argue 
as in \cite{AS2, Lemma 5.4}. Indeed, let us consider the adjoint action 
of $\Gamma$ on $A$. The decomposition $A_1 = A(0) \oplus A(1)$ is 
of $\Gamma$-modules. By \cite{AS2, Lemma 3.1} $A(0)$ is the isotypical 
component of trivial type of the $\Gamma$-module $A_{1}$ (here we use that 
$A$ is finite dimensional). Hence $A(1)$ is the sum of all the  other 
isotypical components. This implies that $\Psi(A(1)) = B(1)$, since 
$\Psi$ intertwines the respective adjoint actions of $\Gamma$. It is 
not difficult to see that $A$ is generated by $A_{1}$ as an algebra, 
by definition of Nichols algebra and \cite{AS2, Lemma 2.4 (iii)}. Hence  
$\Psi(A(n)) = B(n)$ for all $n \ge 0$. We can apply now (the second half 
of) Lemma 6.2.  The automorphism $\psi \in \Aut (H)$ is determined 
by  $\varphi\in \Aut (\Gamma)$. We also know that $R \simeq S^{\psi}$ 
as braided Hopf algebras; but this is equivalent to $V \simeq W^{\psi}$ 
in $_{H}^{H}\Cal{YD}$. This is possible if and only if (6.3) holds.   
If (6.3) holds, we do not need to assume that $R$ and $S$ are 
finite dimensional. By what we have just said, $R \simeq S^{\psi}$. 
Hence $R\# H \simeq S^{\psi}\# H$. By Lemma 6.1, we infer that 
$A\simeq B$.  \qed\enddemo    

\bigpagebreak We now finish the proof of Theorem 1.3. 
We shall apply Proposition 6.3 with $\Gamma = \Bbb Z/(p)$. 
Let first assume $\theta = 2$.   Let $g(1), g(2) \in \Gamma$ and 
$\chi(1), \chi(2) \in \VGamma$  define a Nichols algebra of finite Cartan type.  
Let $u \in \Gamma - 1$, $q \in k^{\times}$ of order $p$ and integers $b$, 
$d$ such that  the $g(i)$'s and $\chi(j)$'s are given by (6.1).  
Note that $d$  is determined by $b$. We denote by $R(q,b)$ the 
Nichols algebra corresponding to the $g(i)$'s and $\chi(j)$'s.  Let now $h(1), 
h(2) \in \Gamma$, $\eta(1), \eta(2) \in \VGamma$) define another Nichols algebra of 
finite Cartan type.  Let also $v \in \Gamma - 1$, $\widetilde q \in k^{\times}$ 
of order $p$, integers $s$, $t$ and 

$$ h(1) = v, \quad h(2) = v^{s}, \quad \langle \eta(1),  
v\rangle = \widetilde q, \quad  \langle \eta(2),  v\rangle = \widetilde q^{t}. $$   

We can assume, by Proposition 6.3, that $u = v$.  Then the Nichols algebra 
corresponding to the $h(i)$'s and $\eta(j)$'s is $R(\widetilde q, s)$.  

\proclaim{Lemma 6.4}The Hopf algebras $A = R(q,b)\# k\Gamma$ and 
$B= R(\widetilde q, s)\# k\Gamma$  are isomorphic if and only if  

\roster \item"(i)" $q =\widetilde q$ and $b \equiv s \mod p$  
when the type is $B_{2}$ or $G_{2}$,  

\item"(ii)" $q =\widetilde q$ and $b \equiv s$ or $\equiv -1 - s \mod p$  
when the type is $A_{2}$ . \endroster \endproclaim  

\demo{Proof}We assume that $A$ and $B$ are isomorphic; it is clear that they 
are of the same type. Let $\varphi$, $\sigma$ be as in Proposition 6.3, 
{\it i.e.} (6.3) holds. If $\sigma = \id$, then $\varphi(u) = u$ and
 hence $q = \widetilde q$, $b = s$. If $\sigma\neq \id$, then 
$\varphi(u) = u^{s}$, $\varphi(u^{b}) = u$, 

$$q =  \langle \chi(1),  
u\rangle = \langle \eta(2),  u^{s}\rangle =\widetilde q^{st}, 
\qquad  q^{d} =  \langle \chi(2),  u\rangle = \langle \eta(1),  u^{s}\rangle =\widetilde q^{s}.$$  

These equalities imply $bs \equiv 1 \mod p$, $dt \equiv 1 \mod p$. Hence 
$st \equiv (bd)^{-1} \mod p$. We know that $(bd)^{-1}\equiv(st)^{-1}\equiv 1, 
2 $ or $3 \mod p$ if the type is respectively $A_{2}$, $B_{2}$ or $G_{2}$. 
This leads immediately to a contradiction unless the type is $A_{2}$. But 
in this last case, $s \equiv d \mod p$ and $t \equiv b \mod p$. Hence 
$q =\widetilde q$ and $b \equiv -1 - s \mod p$. \qed\enddemo   

\bigpagebreak We now pass to the cases $\theta = 3$ and $\theta = 4$.  
We denote by $R_{4}(q, e)$, resp. $R_{3}(q, e)$, the Nichols algebra corresponding to 
the sequence in (6.7), resp. the sequence of the first three terms in (6.7). 
Here $q$ is a third root of 1 and $e =1$ or 2. By Proposition 6.3, arguing 
as for $\theta = 2$, we are reduced to prove:   \proclaim{Lemma 6.5} (i). 
The Hopf algebras  $R_{3}(q, e) \# H$, $q$ a third root of 1, $e = 1$ or 2, 
are not isomorphic to each other.  (ii). The Hopf algebras  $R_{3}(q, e) 
\# H$, $q$ a third root of 1, $e = 1$ or 2, are not isomorphic to each other. 
\endproclaim  

\demo{Proof} (i). All the four elements $g(i)$ are equal for $e=1$ but  
$g(1) \neq g(3)$ for $e = 2$. All the four elements $\chi(j)$ are equal 
for $e=2$ but  $\chi(1) \neq \chi(3)$ for $e = 1$. The statement follows 
from these two observations.  (ii). This is completely analogous and is left 
to the reader. \qed\enddemo  

\bigpagebreak   

\bigpagebreak \subhead \S 7. Pointed Hopf algebras of order $p^4$\endsubhead  

Let $p$ be a prime number. 
We now characterize coradically graded  pointed Hopf
algebras of order $p^4$. This result implies the classification
of all  pointed Hopf algebras of order $p^4$ generated by group-like and skew
primitive elements.  
We shall postpone to a separate paper the statement and proof
of this last result.

\proclaim{Theorem 7.1} Let $G$ be a coradically graded  pointed Hopf
algebras of order $p^4$  with coradical $\k\Gamma$. 
Then $G \simeq R\# \k \Gamma$, where 

\roster 
\item"(7.2)" $R$ is a quantum line if $\Gamma$ has order $p^3$; 

\item"(7.3)" $R$ is a quantum line or plane if $\Gamma$ has order $p^2$;

\item"(7.4)" $R = \goth B(V)$, where $V$ has Dynkin diagram $A_{2}$ 
if $\Gamma$ has order $p$, where $p = 3$ or  $p-1$ is divisible by 3.

\endroster \endproclaim 

\demo{Proof}The case (7.2) follows from Proposition 7.5 below; we can
not apply directly \cite{AS2} since there are non-abelian groups of order 
$p^3$. The case (7.3) follows directly from \cite{AS2}. 

Let us consider the case (7.4). Let $V = P(R)$ as usual. We know that 
$\dim P(R) > 1$. Let $R' = \goth B(V)$; it follows from Theorem 1.3
by a dimension argument
that $V$ is as claimed and $R = R'$, or $V$ is a quantum plane and $R'$
has dimension $p^2$.

We now discard the case when $V$ is a quantum plane and $R'$
has dimension $p^2$ because of Proposition 5.4. \qed \enddemo

\proclaim{Proposition 7.5} Let $\Gamma$ be a finite (non-necessarily abelian)
group. Let $R$ be a braided Hopf algebra in $^{\Gamma}_{\Gamma}{\Cal YD}$
of prime order such that $R_0 = \k 1$.  Then $R$ is a quantum line. \endproclaim 

\demo{Proof} As $P(R)$ is a Yetter-Drinfeld submodule of $R$ and it is non-zero by hypothesis, there exists $g\in \Gamma$ such that $P(R)^g \neq 0$. 
But the Yetter-Drinfeld condition shows that $P(R)^g$ is stable by the action 
of $g$. Hence, there exists $y \in P(R)$, $y \neq 0$, such that 
$$g.y = qy, \quad \delta(y) = g \otimes y.$$
The proof follows now exactly as in the final step of \cite{AS1, Thm. B}.
\qed \enddemo

\remark{Remark} Proposition 7.5 allows the classification of finite dimensional
pointed Hopf algebras whose coradical has index $p$. S. Dascalescu told
us that he also knows this last result. \endremark

\subhead \S 8. Hopf algebras generated in degree one \endsubhead 

\subsubhead  Some examples\endsubsubhead

Let us first discuss examples of braided Hopf algebras $R$ satisfying (1.1), 
(1.2) and (1.3) but not (1.4).  It is enough to show examples of braided Hopf algebras 
$S$ satisfying (1.1), 
(1.2) and (1.4) but not (1.3).

\proclaim{Lemma 8.1} Let $R = \oplus_{n\ge 0} R(n)$ be a graded braided Hopf algebra
in $_{\Gamma}^{\Gamma}\Cal{YD}$, $\Gamma$ a finite abelian group such that   
 $R(0) = \k 1$. Assume that the homogeneous components of $R$
are finite dimensional. Let  $S$ be the graded dual of $R$: $S = \oplus_{n \ge 0} R(n)^*$;
it is a graded braided Hopf algebra
in $_{\Gamma}^{\Gamma}\Cal{YD}$ and  $S(0) = \k 1$.
Then $R$ satisfies (1.3) if and only if $S$ satisfies (1.4). \endproclaim
\demo{Proof} This is implicit in \cite{AG, Example 3.2.8}. 
This was also observed by the referee of \cite{AS2}. \qed \enddemo

We show two easy examples:

\medpagebreak
(8.1). Let $F$ be a field of positive characteristic $p$. Let $S$ be the 
(usual) Hopf algebra $F[x]/\langle x^{p^{2}}\rangle$ with $x\in P(S)$. Then $x^p \in P(S)$.  
Hence $R = S^*$  satisfies (1.1), (1.2) and (1.3) but not (1.4).

\medpagebreak
(8.2).
Let $S = \k[X] = \oplus_{n\ge 0} S(n)$ 
be a polynomial algebra in one variable. We consider $S$ as a
braided Hopf  algebra in $_{H}^{H}\Cal{YD}$, where $H = \k \Gamma$, $\Gamma$
an infinite cyclic group with generator $g$, with action, coaction and comultiplication given by
$$
\delta(X^n)= g^n \otimes X^n, \quad g.X = qX, \quad\Delta(X) = X\otimes 1
+ 1 \otimes X.
$$
Here $q \in \k$ is a root of 1 of order $N$. 
That is, $S$ is a so-called quantum line. Then $S$ satisfies (1.1), 
(1.2) and (1.4) but not (1.3) since $X^N$ is also primitive. Hence the
graded dual $R = S^d = \oplus_{n\ge 0} S(n)^*$ is a braided Hopf algebra
satisfying (1.1),  (1.2) and (1.3) but not (1.4). Let us present explicitly
$R$. It is a vector space with basis $y_n$, $n\ge 0$; the 
structure is given by $y_0 = 1$,
$$\gather
y_ny_m = \bmatrix m+n \\ n\endbmatrix_q y_{n+m}, \quad \Delta(y_n) 
= \sum_{0\le i \le n} y_i \otimes y_{n-i}, \\ \delta(y_n) = g^{-n}
\otimes y_n, \quad g.y_n = q^{-n}y_n.
\endgather$$

\medpagebreak
 We do not know any finite dimensional 
counterexample in characteristic zero.  Note that in both
examples above there is implicit a Frobenius homomorphism. 
In characteristic zero, Lusztig said that "there are no powers of the quantum
Frobenius map" \cite{L4, 8.5, p. 58}; this could be related to our lack of examples.
Indeed our argument below for braided Hopf algebras over an abelian group of exponent
$p$ is a version of Lusztig's claim.

\bigpagebreak
\subsubhead  Proof of Theorem 1.4\endsubsubhead

We now prove Theorem 1.4 (see also  Proposition 5.4 for $\Gamma = \Bbb Z/(p)$, $p>3$).
From now on $\Gamma = \left(\Bbb Z/(p) \right)^s$. 
However, for any finite abelian group $\Gamma$ such that (4.8) holds
the arguments below apply.

By \cite{AS2, Lemma 2.2} it is enough to assume that 
$A$ is coradically graded. By the formula of the multiplication in the biproduct
(see for instance \cite{AS2, Section 2}) it is enough to prove that
a finite dimensional braided Hopf algebra $R$ fulfilling (1.1), 
(1.2) and (1.3) also satisfies (1.4). By Lemma 8.1, we are reduced to prove

\proclaim{Theorem 8.2}Let $S = \oplus_{n\ge 0} S(n)$ be a graded braided Hopf algebra
in $_{\Gamma}^{\Gamma}\Cal{YD}$. We assume that $S$ is  finite dimensional and
 $S(0) = \k 1$. If $S$ is generated as an algebra by $S(1)$ then $P(S) = S(1)$.
\endproclaim

\demo{Proof} Let $V = S(1)$; clearly $V \subseteq P(S)$. 
If $\dim V = 1$, we are done by \cite{AS2, Th. 3.2}. We shall assume that $\dim V > 1$.

It is not difficult to verify the existence of a map $\psi: S \to \goth B(V)$
of braided Hopf algebras in $_{\Gamma}^{\Gamma}\Cal{YD}$ which is the identity on $V$
 (for instance, use Lemma 8.1).
By Corollary 1.2, $V$ is of finite Cartan type.
By Lemma 2.3, we can twist $V$ and assume $b_{ij} = b_{ji}$, 
changing the group if necessary. By Lemma 4.3, 
 we then know that $V$ is of finite FL-type.
Let $x_1, \dots, x_{\theta}$ be a basis of $V$ as in (1.5) and let $(a_{ij})$
be the associated Cartan matrix as in (1.9), (1.10). 

\medpagebreak
Let $U_2^+$ be the quotient of the braided Hopf algebra $T(V)$ by the
ideal generated by the elements $z_{ij} := (\ad_c x_{i})^{1 - a_{ij}} x_{j}$, $i\neq j$.
By Lemma A.1, $z_{ij}$ is primitive in $T(V)$; hence its image in $\goth B(V)$
is 0. Note also that $U_2^+$ is a braided Hopf algebra in $_{\Gamma}^{\Gamma}\Cal{YD}$,
since the defining ideal is a Yetter-Drinfeld submodule of $T(V)$. 
Therefore, we have a commutative diagram of braided Hopf algebras in $_{\Gamma}^{\Gamma}\Cal{YD}$

$$ \CD
T(V) @>>> U_2^+ \\ @VVV @VV\phi V \\ S @>\psi>> \goth B(V).
\endCD \tag 8.3$$

In what follows, we denote by the same symbol $x_i$ the element in any of the algebras
of this diagram.

\medpagebreak
\proclaim{Claim 8.3} The map $\phi$ factorizes through $\psi$: $\phi = \zeta\psi$, where 
$\zeta:  U_2^+ \to S$ is a morphism of braided Hopf algebras in $_{\Gamma}^{\Gamma}\Cal{YD}$. \endproclaim

Let $i\neq j$, $y_{1} = x_{i}$, $y_{2} = z_{ij}$ and $T$ the subalgebra
of $S$ generated by $y_{1}$ and $y_{2}$. Note that 
$$ \delta(y_{2}) = g(i)^{1 - a_{ij}}g(j) \otimes y_{2}, \qquad 
h.(y_{2}) = \chi(i)^{1 - a_{ij}} (h) \chi(j)(h)y_{2}, \quad 
\forall h\in \Gamma.
$$
Assume that $y_2 \neq 0$. Since $T$ is finite dimensional, passing to the graded Hopf algebra associated to its coradical filtration, we conclude from Corollary 1.2 that the matrix
$$\bold d = \pmatrix \langle \chi(i), g(i) \rangle & \langle 
\chi(i)^{1 - a_{ij}}\chi(j), g(i) \rangle
\\ \langle \chi(i), g(i)^{1 - a_{ij}}g(j) \rangle & 
\langle \chi(i)^{1 - a_{ij}}\chi(j), g(i)^{1 - a_{ij}}g(j) \rangle\endpmatrix \tag 8.4$$ 
is of finite Cartan type, say with Cartan matrix $(A_{k\ell})_{1 \le k, \ell \le 2}$. Now
$$
 \langle \chi(i)^{1 - a_{ij}}\chi(j), g(i) \rangle \langle \chi(i), g(i)^{1 - a_{ij}}g(j) \rangle 
=  \langle \chi(i), g(i) \rangle^{2 - a_{ij}}  =  \langle \chi(i), g(i) \rangle^{A_{12}}.
$$
Hence $2 - a_{ij} = 2, 3, 4, 5$ is congruent to $A_{12} =0$,  $-1$, $-2$ or $-3$  
mod the order of $ \langle \chi(i), g(i) \rangle$. 
This is a contradiction since $p > 17 > 7$. The claim is proved.

\medpagebreak
\proclaim{Claim 8.4} The map $\psi$ is an isomorphism. \endproclaim
The algebra $U_2^+$ is the positive part of a quantized  algebra at a root of unity
as in many places in the literature.
We shall follow now the exposition in \cite{AJS}, see also \cite{dCP}. 
Let $R$ be the root system 
corresponding to the Cartan matrix $(a_{ij})$, $R^{+}$ a system of positive
roots, $\Sigma$ the system of simple roots in $R^+$. We can numerate $\Sigma$ by $\{1, \dots, 
\theta\}$. We identify $x_{i}$ with the element $E_{\alpha}$, with respect to this numeration.

\medpagebreak
The algebra $U_2^0$ in \cite{AJS, p. 14}, is nothing but the group algebra of the group
$\Bbb Z^{\theta}$. We consider the Hopf algebra $U_2^{\ge 0} := U_2^0U_2^+$; it has a
PBW basis which is a subset of the basis in \cite{AJS, p. 15, formula (1)}. For each
$\beta \in R^+$, there is a "root vector" $E_{\beta} \in U_2^+$. We consider 

\roster \item"(a)" the ideal $I^+$ of $U^+_2$ generated by all $E_{\beta}^p$, $\beta \in R^+$;
\item"(b)" the subalgebra $\goth{U_2^{\ge 0}}$ of $U_2^{\ge 0}$ generated by $E_{\beta}^p$,
$\beta \in R^+$, and all the $K^{\pm p}_{\alpha}$, $\alpha \in\Sigma$;
\item"(c)" the subalgebra $\goth{U_2^{+}}$ of $U_2^{+}$ generated by $E_{\beta}^p$,
$\beta \in R^+$.
\endroster

Then we have:
\roster \item"(i)" $U^+_2/ I^+ \simeq \goth B(V)$.
\item"(ii)" $\goth{U_2^{\ge 0}}$ is a central Hopf subalgebra of $U_2^{\ge 0}$.
\item"(iii)"   $\goth{U_2^{+}}$ is a central Hopf subalgebra of $U_2^{+}$ 
in $_{ U_2^0}^{ U_2^0}\Cal{YD}$ with trivial braiding.
\endroster

\medpagebreak
The statement (i) follows from the considerations in the last paragraph of p. 15 and
the first paragraph of p. 16 in \cite{AJS}, in view of Theorem 3.1. 

\medpagebreak
The statement (ii) is \cite{dCP, \S 19, Corollary in p. 120}, in presence of the Theorem
four lines below the Corollary.

\medpagebreak 
We proceed with the statement (iii). It is clear by (ii)
that  $\goth{U_2^{+}}$ is  central in $U_2^{+}$.
Now $U_2^0$ is a Hopf subalgebra of $U_2^{\ge 0}$ and its inclusion
admits a retraction $U_2^{\ge 0} \to U_2^0$ sending the (non-trivial) monomials in $E_{\beta}$
to 0. It follows {\it e.g.} from the PBW basis that $U_2^{+}$ is a braided Hopf algebra
in $_{ U_2^0}^{ U_2^0}\Cal{YD}$. It is clear that the algebra, coalgebra and braiding
in this setting are exactly the same as those from (8.3).  By (ii), the action
of $U_2^0$ on $\goth{U_2^{+}}$ is trivial; hence the braiding is trivial. Finally, 
$\goth{U_2^{+}}$ is a  Hopf subalgebra of $U_2^{+}$ since $\goth{U_2^{\ge 0}}$ is a 
 Hopf subalgebra of $U_2^{\ge 0}$.

\bigpagebreak We prove now the Claim 8.4; Theorem 8.2 and Theorem 1.4 follow immediately.
Note that $\goth{U_2^{+}}$ is a connected coalgebra since so is $U_2^{+}$. Let
$S' = \zeta(\goth{U_2^{+}})$; this is a connected braided Hopf subalgebra of $S$. 
If $S' \neq \k$ then there exists $x\in P(S')$, $x\neq 0$. But then the powers of $x$ are 
linearly independent: prove this by induction using the comultiplication and the 
triviality of the braiding. This is a contradiction; hence $S'= \k$ and {\it a fortiori} 
$\zeta(\goth{I^{+}}) = 0$. The Claim follows now from (i).
\qed\enddemo

We finally observe that the proof of Theorem 8.2 can be adapted to finish the
proof of Proposition 5.4.

\proclaim{Lemma 8.5} Let $S$ be a finite dimensional 
braided Hopf algebra in $_{\Gamma}^{\Gamma}\Cal{YD}$, $\Gamma = \Bbb Z/(3)$.
Assume that (1.1), (1.2) and (1.4) hold. Then $S$ is a Nichols algebra.
 \endproclaim

\demo{Proof} The considerations before Claim 8.3 are still valid, using 
Lemma 5.2 instead of Corollary 1.2. Note that by  (1.10) all the $a_{ij}$'s
are 0, -1 or -2. We also observe that Claim 8.4 is also valid. That is, it is 
enough to prove Claim 8.3 for $\Gamma = \Bbb Z/(3)$. 
By Proposition 5.1, we can assume that $\dim S(1)$ is 2 or 3. Indeed, 
$\dim S(1) = 4$ is maximal, and $\dim S(1) = 1$ is treated via \cite{AS2, Thm. 3.2}.

We assume first that $\dim S(1)$ is 2. Then there exists a generator $u$ of 
$\Gamma$, $q \in k^{\times}$ of order $p$ and integers $b$, $d$ such that 

$$ g(1) = u, \quad g(2) = u^{b}, \quad \langle \chi(1),  u\rangle 
= q, \quad  \langle \chi(2),  u\rangle = q^{d}. \tag 5.1$$  

So $\dsize b_{11} = q, \quad b_{22} = q^{bd}$.

There are two possibilities: a quantum plane or $A_{2}$. 

In the first, we have to show that the skew-commutator $z_{12}$ is 0. But 
$$\delta(z_{12}) = u^{1 + b} \otimes z_{12}, \qquad u.z_{12} = q^{1 + d}z_{12}.
$$
Now $u^{1 + b}.z_{12} = q^{(1 + b)(1 + d)}z_{12}$. Since
$b = -d$, the exponent $(1 + b)(1 + d) = 1- b^2$ is 0 mod 3. This is a
contradiction unless $z_{12} = 0$. 

In the second, we have to show that a "quantum Serre relation"-type element 
$z_{12}$ is 0. 
In this situation, $b = d = 1$. But 
$\delta(z_{12}) = u^{2 + 1} \otimes z_{12}$. Again, this is a
contradiction unless $z_{12} = 0$. 

We assume finally that $\dim S(1)$ is 3. In this case, $S(1)$ is the direct
sum $V_1 \oplus V_2$, where say $V_{1}$ is of type $A_{2}$ and
$V_{2}$ is of type $A_{1}$, cf. Proposition 5.1. Applying the considerations
for case 2, we see that $S$ is the twisted tensor product of $\goth B(V_{1})$
and $\goth B(V_{2})$. It follows that $S$ is  $\goth B(V)$, as claimed.
\qed\enddemo

\subhead  Appendix. Primitive elements from quantum Serre relations \endsubhead

\subsubhead  Gaussian binomial coefficients\endsubsubhead

In the polynomial algebra $\Bbb Z [ \bold q]$, we consider the Gaussian, 
or
 $\bold q$-binomial,  coefficients 
 $${\binom ni}_{\bold q} = \frac { {(n)!}_{\bold q}} { {
(n-i)!}_{\bold q} {(i)!}_{\bold q} }, $$ where 
${(n)!}_{\bold q} = {(n)}_{\bold q}  \dots {(2)}_{\bold q}
{(1)}_{\bold q}$, and  ${(n)}_{\bold q} = 1 + {\bold q}
+ \dots + {\bold q}^{n-1},$
 for $n \in \Bbb N$, $0\leq i \leq n$.  We have the identity $${\bold q}^h {\binom
nh}_{\bold q} + {\binom n{h-1}}_{\bold q} = {\binom nh}_{\bold q} +
{\bold q}^{n + 1 - h}{\binom n{h-1}}_{\bold q} ={\binom
{n+1}h}_{\bold q}, 1\leq h \leq n; \tag  A.1$$ This implies that ${\binom 
ni}_{\bold q} \in \Bbb
Z [{\bold q}]$. If $A$ is an associative
algebra over $k$ and $q \in k$, then  ${\binom ni}_q$ denotes the
specialization of ${\binom ni}_{\bold q}$ in $q$.  If $x, y \in A$ are
two elements that $q$-commute, i.e. $ xy = q yx $, then the
quantum binomial formula holds for every $n \in \Bbb N$:$$(x+y)^n
= \sum_{i = 0}^n { \binom ni }_q y^i x^{n-i}. \tag A.2$$
Let us also record
$${(r)}_{\bold q} + \bold q^{r}{(s)}_{\bold q} ={(r + s)}_{\bold q} . \tag A.3$$

\bigpagebreak
There is another version of  "$\bold q$-binomial  coefficients":
in  the Laurent polynomial algebra $\Bbb Z [ \bold q, \bold q^{-1}]$,
we consider 
$${\cbinom ni}_{\bold q} = \frac { {[n]!}_{\bold q}} { {
[n-i]!}_{\bold q} {[i]!}_{\bold q} }, \quad \text {where }
{[n]!}_{\bold q} = {[n]}_{\bold q}  \dots {[2]}_{\bold q}
{[1]}_{\bold q}, \quad \text {and } {[n]}_{\bold q} = \frac{\bold
q^n - \bold q^{-n}}{ \bold q - \bold q^{-1}}.$$ Clearly, 
$$ {(n)}_{\bold q^2}=\bold q^{n-1}{[n]}_{\bold q},$$
and hence 
$${\binom ni}_{\bold q^2} = \bold q^{i(n-i)}{\cbinom ni}_{\bold q}.$$
In particular, the following equality- see \cite{L3, 1.3.4 (a)}-

$$ \sum_{i = 0} ^{n} (-1)^i \bold q^{i(1-n)}{\cbinom ni}_{\bold q} = 0 
\tag A.4 $$
translates into
$$ \sum_{i = 0} ^{n} (-1)^i \bold q^{(i^2 + i) / 2 - ni}{\binom 
ni}_{\bold q} = 0. \tag A.5 $$
This last equality is also equivalent, by a change of the summation 
index, to 
$$ \sum_{h = 0} ^{n} (-1)^h \bold q^{(h^2 - h) / 2}{\binom 
nh}_{\bold q} = 0. \tag A.6 $$

 \bigpagebreak
\subsubhead  Primitive elements\endsubsubhead

Let $R$ be a braided Hopf algebra over our fixed  finite abelian group
$\Gamma$. We assume that (1.1), (1.2) and (1.3) hold; we also suppose that
$\dim P(R) = 2$. But we do not assume that $R$ is finite dimensional,
unless explicitly stated. To simplify the notation, we fix a basis $x,y$ of
$P(R)$ such that $$\delta (x) = g\otimes x, \quad h. x = \chi(h)x \quad,
\delta (y) = t\otimes y, \quad h. y = \eta(h)y \quad \forall h\in \Gamma,$$
for some $g,t \in \Gamma$, $\chi, \eta \in \widehat\Gamma$. Let $N$, $M$
be respectively the orders of $\chi(g)$, $\eta(t)$. We do suppose that they
are positive.

\definition{Definition} Let $S$ be a braided Hopf algebra. The {\it
braided adjoint representation\/} is the linear map $\ad_{c}: S\to
\End S$ given by
$$\ad_{c}(u)(v) = \mu(u\otimes v - c(u\otimes v));$$
here $\mu$ is the multiplication map of $S$ and $c$ is the
commutativity constraint. \enddefinition
The following property holds:
$$\ad_{c}(u)(vw) = \ad_{c}(u)(v)w + (u_{(-1)}.v) \ad_{c}(u_{(0)})(w).
\tag A.7$$

\bigpagebreak
In our situation, we define inductively  $z_{1} =\ad_{c}(x)(y)= xy -
\eta(g)yx$;  $z_{j + 1} = \ad_{c}(x)(z_{j})$. For instance, 
$$\align z_{2} &= x^{2}y - \eta(g)\left(1+\chi(g)\right) xyx +
\eta(g)^{2}\chi(g) yx^{2}, \\  
z_{3} &= x^{3}y - \eta(g) (3)_{\chi(g)} x^{2}yx +  \eta(g)^{2}\chi(g)
(3)_{\chi(g)} xyx ^{2} -    \eta(g)^{3}\chi(g)^{3}yx ^{3}\endalign$$
and in general
$$z_{N} = \sum_{i=0}^{N} (-1)^{i}  {\binom Ni}_{\chi(g)}
\chi(g)^{i(i-1)/2} \eta(g)^{i}\, x^{N - i}yx^{i}.  \tag A.8$$
The proof of (A.8) is by induction using (A.7).

We following Lemma is probably well-known. For completeness, we give a direct proof.
\proclaim{Lemma	A.1}  Assume that
$$  \eta(g) \chi(t) \chi(g)^{r-1} = 1. \tag A.9$$
 Then $(\ad_c x)^r y $ is primitive.
\endproclaim 

\demo{Proof}Let $z= \sum_{i=0}^r\alpha_i x^i y x^{r-i}$, where $\alpha_i 
\in k$. Then 
$$\align
\Delta(z) &= \sum_{i=0}^r \sum_{\ell=0}^i\sum_{h=0}^{r-i} 
\alpha_i {\binom i \ell}_{\chi(g)}{\binom {r-i}h}_{\chi(g)}
\left[\eta(g)^{i - \ell}\chi(g)^{h(i - \ell)}  x^{\ell} y x^{h}  \otimes 
x^{r - \ell - h}   \right. \\
&+ \left. \chi(t)^{h}\chi(g)^{h(i - \ell)}  x^{\ell + h} \otimes x^{i -\ell} y 
x^{r - i- h} \right] \\
&= z \otimes 1 + 1 \otimes z  \\ &  + \sum \Sb \ell \ge 0, h\ge 0 \\ \ell + h < r \endSb \left( \sum_{\ell \le i \le r -h} \alpha_i {\binom i \ell}
{\binom {r-i}h} \eta(g)^{i - \ell}\chi(g)^{h(i - \ell)}  \right) x^{\ell} y 
x^{h}  \otimes x^{r - \ell - h}   \\ &   + 
 \sum \Sb u \ge 0, v\ge 0 \\ u + v < r \endSb \left( \sum_{u \le i \le r- v} 
\alpha_i {\binom i u}{\binom {r-i}v} \chi(t)^{r - i - v}
\chi(g)^{u(r - i - v)}  \right) x^{r - u - v}\otimes x^{u} y x^{v}.   
\endalign$$

Hence, if $\alpha_0, \dots, \alpha_r$ is a solution of the system of
equations
$$\align
\sum_{\ell \le i \le r -h} \alpha_i {\binom i \ell}
{\binom {r-i}h} \eta(g)^{i - \ell}\chi(g)^{h(i - \ell)} &= 0, \tag A.10a\\
\sum_{u \le i \le r- v} 
\alpha_i {\binom i u}{\binom {r-i}v} \chi(t)^{r - i - v}
\chi(g)^{u(r - i - v)} &= 0, \tag A.10b
\endalign$$
then $z$ is primitive.

\bigpagebreak
Now assume that $z = (\ad_c x)^r y$, that is 
$$ z = \sum_{0 \le i \le r} (-1)^{r - i}  {\binom  ri}_{\chi(g)}
\chi(g)^{(r- i) (r - i - 1) / 2} \eta(g)^{ r- i} x^iyx^{r-i}. $$
To see that $z$ is primitive, we  divide $z$ by
$(-1)^{r}\chi(g)^{(r^2  - r) / 2} \eta(g)^{r}$ and reduce to check that 
$$\alpha_i :=  (-1)^{i}  {\binom  ri}_{\chi(g)}
\chi(g)^{(i^2 + i)/ 2 - ri} \eta(g)^{-i}$$
satisfies the system of equations (A.10). We proceed with (A.10a). 
If $\ell \ge 0$, $h \ge 0$, $\ell + h < r$, then

$$ \multline 
\sum_{\ell \le i \le r-h} (-1)^{i} {\binom  ri}_{\chi(g)}
 {\binom i \ell}_{\chi(g)}{\binom {r-i}h}_{\chi(g)} 
\eta(g)^{i - \ell - i}\chi(g)^{(i^2 + i) / 2 - ri + h(i - \ell)}
\\ =  \frac{(r)! \eta(g)^{-\ell}\chi(g)^{-h\ell}}{(\ell)! (h)! (r - \ell -h )!}
\sum_{\ell \le i \le r-h} (-1)^{i}  \frac{(r - \ell -h)!}{(i-\ell)!(r-i-h)!}
\chi(g)^{(i^2 + i) / 2 + i(h - r)}.
\endmultline$$
Dividing by the appropiate factor and changing the summation index, we get
$$(-1)^{\ell} \chi(g)^{(\ell^2 + \ell) / 2}  \sum_{0 \le i \le r - h - \ell}
 {\binom {r - h - \ell} i}_{\chi(g)} \chi(g)^{(i^2 + i) / 2 
+ i(h + \ell - r)} = 0, $$
by (A.5).

Now we pass to (A.10b). For $u \ge 0$, $vv \ge 0$, $u + v < r$, we compute:
$$ \multline 
\sum_{u \le i \le r-v} (-1)^{i} {\binom  ri}_{\!\!\chi(g)}
 {\binom i u}_{\!\!\chi(g)}{\binom {r-i}v}_{\!\!\chi(g)} 
\eta(g)^{-i}\chi(t)^{r - i - v}
\chi(g)^{(i^2+i) / 2 - ri + u(r - i - v)}\mkern-10mu
\\
 = \frac{(r)! \chi(t)^{r-v}    \chi(g)^{u(r-v)}}{(u)! (v)! (r - u -v )!}
\!\! \sum_{u \le i \le r-v}\!\! (-1)^{i}  \frac{(r - u -v)!}{(i-u)!(r-i-v)!}
\chi(g)^{(i^2+i) / 2  -ri -ui +i(r-1)},\mkern-10mu
\endmultline$$
where we have used (A.6). After dividing by the appropiate factor
and changing the summation index, we arrive to
$$
(-1)^{u}\chi(g)^{(u^2  + u) / 2} \sum_{0 \le i \le r-v-u} (-1)^{i}  {\binom {r -
u - v} i}_{\chi(g)} \chi(g)^{(i^2 - i) / 2   } = 0, $$
by (A.6).  \qed\enddemo

   \enddocument